\documentclass[10pt,a4paper]{article}

\usepackage{amssymb,amsmath}
\usepackage{amsthm}

\usepackage{xcolor}
\usepackage{tikz}
\usetikzlibrary{fadings}
\usetikzlibrary{patterns}
\usetikzlibrary{calc}

\newtheorem{thm}{Theorem}[section]

\newtheorem{prop}[thm]{Proposition}
\newtheorem{cor}[thm]{Corollary}
\newtheorem{defn}[thm]{Definition}
\theoremstyle{definition}
\newtheorem{example}[thm]{Example}
\theoremstyle{remark}
\newtheorem{rem}[thm]{Remark}

\newcommand{\Skindef}{[\,\cdot,\cdot\,]}

\newcommand{\cA}{\mathcal A}
\newcommand{\cD}{\mathcal D}
\newcommand{\cG}{\mathcal G}
\newcommand{\cH}{\mathcal H}
\newcommand{\cK}{\mathcal K}
\newcommand{\cL}{\mathcal L}
\newcommand{\CC}{\mathbb C}
\newcommand{\NN}{\mathbb N}
\newcommand{\RR}{\mathbb R}
\newcommand{\fra}{\mathfrak{a}}

\newcommand{\frd}{\mathfrak{d}}

\newcommand{\frt}{\mathfrak{t}}
\newcommand{\defeq}{\mathrel{\mathop:}=}
\newcommand{\defequ}{\mathrel{\mathop:}\hspace*{-0.72ex}&=}
\newcommand{\sess}{\sigma_{\rm ess}}

\newcommand{\amin}{a_{\textup{min}}}
\newcommand{\deltaess}{\delta_{\textup{ess}}}

\DeclareMathOperator{\Real}{Re}

\renewcommand{\Re}{\Real}

\DeclareMathOperator{\ran}{ran}

\DeclareMathOperator*{\essinf}{ess\,inf}
\DeclareMathOperator*{\esssup}{ess\,sup}
\DeclareMathOperator{\essran}{ess\,ran}

\newcommand{\dom}{{\cal D}}

\newcommand{\hhalf}{H_{\frac{1}{2}}}
\newcommand{\half}{\frac{1}{2}}
\newcommand{\mphhalf}{H_{-\frac{1}{2}}\times H_{\frac{1}{2}}}
\newcommand{\eps}{\varepsilon}

\newcommand{\wt}{\widetilde}
\newcommand{\vect}[2]{\begin{pmatrix} #1 \\ #2 \end{pmatrix}}
\newcommand{\smvect}[2]{\bigl(\begin{smallmatrix} #1 \\ #2 \end{smallmatrix}\bigr)}

\newcommand{\be}{\begin{equation}}
\newcommand{\ee}{\end{equation}}
\newcommand{\bea}{\begin{eqnarray}}
\newcommand{\eea}{\end{eqnarray}}
\newcommand{\beann}{\begin{eqnarray*}}
\newcommand{\eeann}{\end{eqnarray*}}

\newcommand{\myassumptionbf}[2]{
\begin{list}{}{\setlength{\leftmargin}{8ex}
\setlength{\labelwidth}{8ex}
\setlength{\labelsep}{2ex}
}
\item[\textbf{(#1)}]
#2
\end{list}
}

\numberwithin{equation}{section}

\title{Variational principles for self-adjoint operator functions arising from
second order systems: estimates below the essential spectrum}%
\author{Birgit Jacob, Matthias Langer, and Carsten Trunk}%

\begin{document}
\maketitle
\begin{abstract}
We investigate quadratic operator functions associated with damped second-order evolution equations in Hilbert spaces. 
Under suitable assumptions on the stiffness and damping operators, we study the relationship between the 
quadratic operator function and the corresponding first-order operator matrix, which is self-adjoint in an associated Krein space.
Our main objective is to establish variational characterizations of real eigenvalues lying below the essential spectrum. 
Using generalized Rayleigh functionals and the structure of the associated operator matrix, 
we derive min-max and max-min principles for two classes of eigenvalues. 
These results are then applied to obtain comparison theorems and explicit eigenvalue bounds for damped second-order systems.

\end{abstract}


\section{Introduction}

Second--order differential equations of the form
\begin{equation}\label{sys}
  \ddot z(t)+D\dot z(t)+A_0z(t)=0
\end{equation}
arise in many applications in mechanics and engineering. Such equations are frequently
used to model vibrating structures in the presence of damping, for instance
transverse oscillations of beams or plates. A typical example is the
Euler--Bernoulli beam with Kelvin--Voigt damping (see, e.g.\ \cite{BSY,CLL,JTW,Pu}),
which can be described by a partial differential equation of the form
\begin{equation}
  \frac{\partial^2u}{\partial t^2}
  +\frac{\partial^2}{\partial r^2}\left[
      \frac{\partial^2u}{\partial r^2}
      +C_d(r)\frac{\partial^3u}{\partial r^2\partial t}
    \right]=0,
 \qquad r\in(0,1),\ t>0.
 \label{eq:beam}
\end{equation}
Here $u(r,t)$ denotes the deflection of the beam from
its rigid body motion at time $t$ and position $r$. Moreover,
$C_d\in L^\infty(0,1)$ with $C_d(r)\ge c>0$ a.e.\ describes the
damping properties of the piezoelectric film.
Appropriate boundary conditions at $r=0$ and $r=1$ complete the beam model.

The Euler--Bernoulli beam \eqref{eq:beam} can be treated as a second-order
problem in the Hilbert space $H=L^2(0,1)$ of the abstract form displayed above.
Here $A_0$ stands for the stiffness and
the Kelvin--Voigt damping is represented by a differential operator rather than a multiplication operator.
Formally, the stiffness and damping operators are given by
\[
  A_0u=\frac{\partial^4u}{\partial r^4},
  \qquad
  Du=\frac{\partial^2}{\partial r^2}
      \left(C_d(r)\frac{\partial^2u}{\partial r^2}\right),
\]
with their precise realizations determined by the boundary conditions.

In this paper we study equations of this type in an abstract Hilbert space
$H$. We assume that the stiffness operator $A_0$ is a self-adjoint,
uniformly positive operator and that the damping operator $D$ satisfies
that $A_0^{-1/2}DA_0^{-1/2}$ is bounded and non-negative in $H$.
More specifically, we assume that $A_0^{-1}$ is compact and that
$A_0^{-1/2}DA_0^{-1/2}$ is uniformly positive (cf.\ the  standing
assumptions (A1)--(A4) below). These hypotheses allow the essential spectrum of the
associated first-order operator to be described in terms of the normalized
damping operator.
This class of systems includes many models arising from partial differential
equations describing damped vibrations.

A common approach for analysing second-order systems is to rewrite them as
first-order evolution equations. Introducing the state vector
\[
  x(t)=\begin{pmatrix}z(t)\\[0.5ex] \dot z(t)\end{pmatrix},
\]
the equation can be transformed into the first-order system
\[
  \dot x(t)=\mathcal A x(t),
\]
where
\[
  \mathcal A=\begin{pmatrix}0&I\\[0.5ex]-A_0&-D\end{pmatrix}
\]
is an operator matrix acting in an appropriate state space.
The operator $\mathcal A$ is the generator of a strongly continuous
semigroup describing the evolution of the system, and its spectral
properties determine essential features of the dynamics such as stability
and decay behaviour of solutions; see, e.g., \cite{jatr,JT08,JTW}.
Although $\mathcal A$ is not self-adjoint in the energy Hilbert space,
it is self-adjoint with respect to a suitable indefinite inner product.
This Krein-space structure is one of the ingredients used to obtain
variational information about its real eigenvalues.

The operator matrix associated with damped second-order systems has been
extensively studied in the literature. Such models occur in stabilization
problems, Riccati equations, and control theory; see, for example, \cite{AI,JTW}.
A fundamental problem is the investigation of the spectrum of the operator
$\mathcal A$ and its relation to the properties of the generated semigroup.
In particular, analyticity of the semigroup and the existence of Riesz
bases of eigenvectors and associated vectors play an important role in
the analysis of these systems \cite{JTW}.

Closely related to the operator matrix formulation is the quadratic
operator function
\begin{equation}\label{Baby}
  T(\lambda)=\lambda^2 I+\lambda D+A_0.
\end{equation}
This operator polynomial arises when replacing the time derivative in the
second-order equation by a complex parameter $\lambda$.
On the regions considered below, the spectral problems for the operator matrix $\mathcal A$ 
and the quadratic operator function $T$ can be identified.
Hence the spectral analysis of the operator matrix $\mathcal A$ can be
studied via the associated spectral problem for $T$;
see \cite{JLT16}.

Variational principles provide a powerful tool for studying eigenvalues of
operators and operator functions. For self-adjoint operators, eigenvalues
below the essential spectrum can be characterized using the classical
Rayleigh functional
\[
  p(x)=\frac{\langle  A x,x\rangle}{\langle x,x\rangle},
\]
together with minimax principles introduced by Rayleigh, Ritz, Courant and
others; see \cite{CourantHilbert}. These techniques play an important role both in qualitative
spectral analysis and in numerical approximation of eigenvalues.

For operator functions arising from second-order systems, generalized
Rayleigh functionals can be introduced in a similar way. Variational
principles for quadratic operator polynomials were first established by
Duffin \cite{Duffin1955} and later extended to more general classes of operator
functions; see, for example, \cite{JLT16}. These results allow the
characterization of eigenvalues of operator functions by minimax principles
analogous to those known for self-adjoint operators.
The generalized Rayleigh functionals $p_+$ and $p_-$ are the two real roots,
whenever they exist, of
\[
  \mathfrak t(\lambda)[x]
  \defeq \lambda^2\|x\|^2+\lambda\mathfrak d[x]+\mathfrak a_0[x]=0.
\]
The distinction between these two branches is important for describing
real eigenvalues in different spectral regions.

The earlier paper \cite{JLT16} established variational characterizations for
real eigenvalues of the first-order operator matrix lying \emph{above} its
essential spectrum, in a spectral region near the origin. Those results
primarily use the upper branch $p_+$ of the generalized Rayleigh
functional and also lead to eigenvalue comparison theorems. The general
min--max framework for parameter-dependent self-adjoint operator functions,
including increasing and decreasing crossings, was already recalled and
adapted to second-order systems in \cite{JLT16}, building on earlier work. 
The present article therefore does not aim to introduce a new
abstract variational principle; instead, it applies and develops this
framework for a different spectral region.

Our main objective is to investigate real eigenvalues \emph{below},
that is, to the left of, the essential spectrum of $\mathcal A$. In
contrast to the spectral region studied in \cite{JLT16}, the present analysis
requires a uniform positivity assumption on the normalized damping
operator. In particular, if
\[
  \delta_{\mathrm{ess}}
  \defeq\min\sigma_{\mathrm{ess}}
  (A_0^{-1/2}DA_0^{-1/2})>0,
\]
then the left endpoint of the essential spectrum of $\mathcal A$
is $-1/\delta_{\mathrm{ess}}$. Theorem~4.2 provides min--max and
max--min characterizations for two branches of real eigenvalues lying
to the left of this endpoint, using $p_+$ and $p_-$, respectively.
The formulae involve index shifts determined by the number of positive
squares of $T(c)$. Moreover, an infinite sequence of real eigenvalues
tends to $-\infty$, while, under suitable spectral conditions, a second
sequence may accumulate at the left endpoint of the essential spectrum.

We then extend the eigenvalue-comparison approach of \cite{JLT16} to these
two branches. Theorem~5.1 describes their behaviour under an ordering of
the stiffness and damping forms. Theorem~5.3 yields explicit two-sided
bounds for the eigenvalues in terms of the eigenvalues of $A_0$ and
spectral bounds on $A_0^{-1/2}DA_0^{-1/2}$. The proportional damping
model $D=dA_0$ illustrates the two distinct asymptotic behaviours:
one sequence tends to $-\infty$, whereas another approaches the left
endpoint of the essential spectrum.
These results allow a direct application for the Kelvin--Voigt damping.

The paper is organized as follows. In Section~\ref{framework} we introduce the
standing assumptions and investigate the quadratic operator function and
its associated sesquilinear forms. Section~\ref{sec:matrix} is devoted to the corresponding
operator matrix and its spectral properties in the associated Krein space.
In Section~\ref{sec:varprinciple} we establish variational principles for real eigenvalues below
the essential spectrum. Finally, Section~\ref{sec:comparison} contains comparison theorems
and explicit eigenvalue estimates, and
the Kelvin--Voigt damping is investigated in Section~\ref{sec:beam}.

\paragraph{Notation.}
By $\mathbb N$ we denote the set of natural numbers,
$\mathbb N=\{1,2,3,\ldots\}$.
A closed, densely defined operator $S$ in a Hilbert space is called Fredholm
if the dimension of its kernel and the (algebraic) codimension of its range
are finite. Denote by $\sigma_p(S)$ the point spectrum, by $\sigma(S)$
the spectrum, and by $\sigma_{\mathrm{ess}}(S)$ the essential spectrum of
$S$,
\[
 \sigma_{\mathrm{ess}}(S)
  \defeq\{\lambda\in\mathbb C:S-\lambda I\text{ is not Fredholm}\}.
\]
The notion of sectorial and m-sectorial operators will play an essential
role; for a definition we refer to \cite[\S V.3.10]{K}. Related to sectorial
operators is the notion of
sectorial forms; see \cite[\S VI.1.2]{K}.
We make frequent use of spectral subspaces $L_I(A)$ of a self-adjoint
operator $A$ in a Hilbert space corresponding to a real interval $I$.
\medskip

\section{The operator polynomial
$\lambda^2 + \lambda D + A_0$ and the related form}\label{framework}

In this section we present the standing assumptions \textbf{(A1)--(A4)} and discuss first consequences.
We assume throughout this paper the following;
compare also with \cite{WT,TWII} and \cite{JLT16,jatr,JT08,JTW}.

\myassumptionbf{A1}{
The stiffness operator $A_0 : \dom (A_0) \subset H \rightarrow H$
is a self-adjoint, positive definite linear operator on an infinite-dimensional Hilbert space $(H, \langle \cdot, \cdot\rangle)$ such that $0\in\rho(A_0)$.
}

\noindent
Here $\langle \cdot, \cdot\rangle$ is the inner product in $H$ and $\|\cdot \|$ is the induced norm. We define
\begin{equation*}
  H_{\frac{1}{2}} \defeq \dom (A_0^{1/2})\text{ with the norm }
  \|x\|_{H_{\frac{1}{2}}} \defeq \|A_0^{1/2} x\|
\end{equation*}
and
\begin{equation*}
  \begin{aligned}
    & H_{-\frac{1}{2}} \text{ as the completion of $H$ with respect to} \;\; \|x\|_{H_{-\frac{1}{2}}} \defeq \|A_0^{-1/2}x\|.
  \end{aligned}
\end{equation*}
The corresponding inner products are denoted by 
$\langle \cdot , \cdot\rangle_{H_{\frac{1}{2}}}$
and $\langle \cdot , \cdot\rangle_{H_{-\frac{1}{2}}}$, respectively. By continuity, $A_0$ and $A_0^{1/2}$ can be extended to isometric
isomorphisms from $H_{\frac{1}{2}}$ onto $H_{-\frac{1}{2}}$ and from $H$
onto $H_{-\frac{1}{2}}$, respectively.  To keep notation simple, these extensions are also denoted
by $A_0$ and $A_0^{1/2}$.  The space $H_{-\frac{1}{2}}$ can
be identified with the dual space of $H_{\frac{1}{2}}$ by identifying elements
$x\in H_{-\frac{1}{2}}$ with functionals on $H_{\frac{1}{2}}$ as follows
\[
  \langle x,y\rangle_{\mphhalf} \defeq \bigl\langle A_0^{-1/2}x,A_0^{1/2}y\bigr\rangle,
  \qquad y\in H_{\frac{1}{2}}.
\]
Note that, for $x\in H$, $y\in H_{\frac{1}{2}}$, we have
\begin{equation*}
  \langle x,y\rangle_{H_{-\frac{1}{2}}\times H_{\frac{1}{2}}}=\langle x,y\rangle.
\end{equation*}
\myassumptionbf{A2}{
The damping operator $D: H_{\frac{1}{2}} \rightarrow
H_{-\frac{1}{2}}$ is a bounded operator with
\[
  \langle Dz, z\rangle_{H_{-\frac{1}{2}}\times
  H_{\frac{1}{2}}} \ge 0  ,\qquad z\in H_{\frac{1}{2}}.
\]
}
\noindent
The condition \textbf{(A2)} is equivalent to the fact that the operator
$D: H_{\frac{1}{2}} \rightarrow H_{-\frac{1}{2}}$ is such that
$A_0^{-1/2}D A_0^{-1/2}$ is a bounded self-adjoint non-negative operator in~$H$. 
This implies that the essential spectrum of the latter operator is non-empty and bounded.  Set
\begin{alignat}{2}\label{Elbersfeld}
  \delta &\defeq \min\sigma\bigl(A_0^{-1/2}DA_0^{-1/2}\bigr), &
  \gamma &\defeq \max\sigma\bigl(A_0^{-1/2}DA_0^{-1/2}\bigr),
  \\[1ex]\label{Elbersfeld2}
  \delta_{\rm ess} &\defeq \min \sess\bigl(A_0^{-1/2}DA_0^{-1/2}\bigr), \qquad &
  \gamma_{\rm ess} &\defeq \max \sess\bigl(A_0^{-1/2}DA_0^{-1/2}\bigr).
\end{alignat}

Note also that, for $x,y\in H_{\frac{1}{2}}$,
\begin{align}
  \langle Dx,y\rangle_{H_{-\frac{1}{2}}\times H_{\frac{1}{2}}}
  &= \bigl\langle A_0^{-1/2}DA_0^{-1/2}A_0^{1/2}x,A_0^{1/2}y\bigr\rangle
  \notag\\[1ex]
  &= \bigl\langle A_0^{1/2}x,A_0^{-1/2}DA_0^{-1/2}A_0^{1/2}y\bigr\rangle
  = \overline{\langle Dy,x\rangle_{H_{-\frac{1}{2}}\times H_{\frac{1}{2}}}}\,.
  \label{D_symm}
\end{align}

\myassumptionbf{A3}{
The operator $A_0^{-1}$ is a compact operator in $H$.
}

\myassumptionbf{A4}{We have $0\in\rho\bigl(A_0^{-1/2}DA_0^{-1/2}\bigr)$.}

\noindent
The numbers defined in \eqref{Elbersfeld}, \eqref{Elbersfeld2} satisfy
\begin{equation*}
  0 < \delta \le \delta_{\rm ess} \le \gamma_{\rm ess} \le \gamma<\infty.
\end{equation*}

\medskip
Within the setting  \textbf{(A1)}--\textbf{(A4)} the operator family $T$ in \eqref{Baby} acts between different spaces, namely between $H_{\frac{1}{2}}$ and $H_{-\frac{1}{2}}$. In order to obtain a setting within $H$ we introduce two  sesquilinear forms $\fra_0$ and $\frd$ in $H$ with
domains $\cD(\fra_0)=\cD(\frd)=H_{\frac{1}{2}}$, defined via
\begin{equation}\label{Wuppi34}
  \begin{aligned}
    \fra_0[x,y] \defequ \langle A_0x,y\rangle_{H_{-\frac{1}{2}}\times
    H_{\frac{1}{2}}}, \\[1ex]
    \frd[x,y] \defequ \langle Dx,y\rangle_{H_{-\frac{1}{2}}\times
    H_{\frac{1}{2}}},
  \end{aligned}
  \qquad x,y\in H_{\frac12},
\end{equation}
and the corresponding quadratic forms are denoted by $\fra_0[x]$ and $\frd[x]$, respectively. 
By \textbf{(A1)}, we have  
\[
  \fra_0[x,y] = \langle A_0^{1/2}x, A_0^{1/2}y\rangle
\]
and, hence,
(see, e.g.\ \cite[Example~VI.1.13]{K})
the form $ \fra_0$ is a closed, symmetric, non-negative
form in $H$. By  \textbf{(A2)},  the form $ \frd$ is  symmetric (see \eqref{D_symm}) 
and non-negative, and, obviously, both forms are densely defined in $H$. 
Moreover, for $x\in H_\half$ we have
\begin{align}
  \frd[x] &= \langle Dx,x\rangle_{H_{-\frac{1}{2}}\times H_{\frac{1}{2}}}
  = \bigl\langle A_0^{-1/2}DA_0^{-1/2}A_0^{1/2}x,A_0^{1/2}x\bigr\rangle 
  \nonumber\\[1ex]
  &\le \gamma\bigl\langle A_0^{1/2}x,A_0^{1/2}x\bigr\rangle
  = \gamma\fra_0[x]
  \label{nuevo}
\end{align}
and similarly,
\begin{equation}\label{nuevo2}
  \frd[x] \ge \delta\fra_0[x].
\end{equation}

For $\lambda\in\CC$ we define the sesquilinear form $\frt(\lambda)$ with
domain $\cD(\frt(\lambda))=H_{\frac{1}{2}}$ by
\begin{align}
  \frt(\lambda)[x,y] \defeq&\; \langle \lambda^2 x+\lambda Dx+A_0x,y
  \rangle_{H_{-\frac12}\times H_{\frac12}}
  \nonumber\\
  =&\; \lambda^2\langle x,y\rangle
  + \lambda\langle Dx,y\rangle_{H_{-\frac12}\times H_{\frac12}}
  + \langle A_0x,y\rangle_{H_{-\frac12}\times H_{\frac12}}
  \nonumber\\
  =&\; \lambda^2\langle x,y\rangle + \lambda \frd[x,y]
  + \fra_0[x,y],
  \hspace*{17ex} x,y\in H_{\frac12},
  \label{Wuppi5}
\end{align}
and the corresponding quadratic form by $\frt(\lambda)[x] \defeq \frt(\lambda)[x,x]$
for $x\in H_{\frac12}$.  Set
\begin{equation}\label{defOmegadelta0}
  \Omega \defeq 
    \biggl\{z\in\mathbb C\,\Bigm|\, \Re z < -\frac{1}{\delta_{\textup{ess}}}\biggr\}.
\end{equation}

In the following theorem we collect some properties of the 
forms $\fra_0$, $\frd$ and $\frt$ defined in \eqref{Wuppi34} and \eqref{Wuppi5} for $\lambda\in\Omega$.
 
\begin{thm}\label{rel_bd}
Let $\lambda\in\Omega$. Then the form $-\lambda\frd$
is closed and sectorial.
Moreover, $\fra_0$ is relatively bounded with respect to $-\Re\lambda\,\frd$ with
$(-\Re\lambda\,\frd)$-bound less than $1$, that is,
there exist constants $a\ge0$, $0<b<1$ 
such that,
\begin{equation}\label{a0_rel_bdd_d}
  \fra_0[x] \le a\|x\|^2 - b \Re\lambda\,\frd[x],
  \quad \mbox{for} \quad x\in H_{\frac{1}{2}}.
\end{equation}

The form $-\frt(\lambda)$
is a closed sectorial form in $H$, which is densely defined.
\end{thm}

\begin{proof}
In order to prove the closedness of $\frd$,
let $(x_n)_{n\in \mathbb N}$ be a sequence in $\dom(\frd)$ which
converges in $H$ to some $x$ such that $\frd[x_n-x_m] \to 0$ for $n,m
\to \infty$.
Then, by \eqref{nuevo2}, $(x_n)_{n\in \mathbb N}$
is a Cauchy sequence in $ H_{\frac{1}{2}}$; hence it converges in
$ H_{\frac{1}{2}}$ to $x$, that is, $x\in  H_{\frac{1}{2}}=\dom(\frd) $.
Moreover, \eqref{nuevo} implies that $\frd[x_n-x] \to 0$ as $n\to \infty$.

Let $\lambda\in\Omega$. We
show that the form $\fra_0$ is relatively bounded with respect
to $-\Re\lambda\,\frd$. Choose $\delta'$ with
$-\frac{1}{\Re \lambda} < \delta'< \delta_{\rm ess}$ and set $\Delta_1 \defeq [0,\delta']$.
Since $\delta'<\delta_{\textup{ess}}$, the spectral subspace
\[
  \cL_{\Delta_1}(A_0^{-1/2}DA_0^{-1/2})
\]
is finite-dimensional and the corresponding spectral projection $E(\Delta_1)$ is of finite rank.
Choose $\eps>0$ such that
\[
  \eps < \Bigl(-\Re\lambda-\frac{1}{\delta'}\Bigr)\delta.
\]
By \cite[Lemma~4.1]{JLT16} there exists a constant $C\ge0$ such that,
for $x\in H_{\frac{1}{2}}$,
\begin{equation*}
  \bigl\langle E(\Delta_1) A_0^{1/2}x,A_0^{1/2}x\bigr\rangle
  = \bigl\|E(\Delta_1)A_0^{1/2}x\bigr\|^2
  \le C\|x\|^2 + \eps \bigl\|A_0^{1/2}x\bigr\|^2.
\end{equation*}
Using \eqref{nuevo2} we obtain, for $x\in H_{\frac{1}{2}}$, that
\begin{align*}
  \fra_0[x]
  &= \langle A_0x,x\rangle_{H_{-\frac{1}{2}}\times H_{\frac{1}{2}}} \\[1ex]
  &= \bigl\langle (I-E(\Delta_1))A_0^{1/2}x,A_0^{1/2}x\bigr\rangle
    + \bigl\langle E(\Delta_1) A_0^{1/2}x,A_0^{1/2}x\bigr\rangle
  \\[1ex]
  &\le \frac{1}{\delta'} \bigl\langle A_0^{-1/2}DA_0^{-1/2}A_0^{1/2}x,A_0^{1/2}x\bigr\rangle
    + C\|x\|^2 + \eps\bigl\|A_0^{1/2}x\bigr\|^2
  \\[1ex]
  &\le \frac{1}{\delta'}\frd[x] + C\|x\|^2 + \frac{\eps}{\delta}\frd[x]
  = \frac{1}{-\Re\lambda}\Bigl(\frac{1}{\delta'}+\frac{\eps}{\delta}\Bigr)(-\Re \lambda)\frd[x] + C\|x\|^2.
\end{align*}
From the definition of $\eps$ we see that
\[
  -\frac{1}{\Re \lambda}\Bigl(\frac{1}{\delta'}+\frac{\eps}{\delta}\Bigr) < 1,
\]
which shows that the form $\fra_0$ is relatively bounded with respect
to $-\Re\lambda\,\frd$ with $(-\Re\lambda\,\frd)$-bound less than $1$.
By  \cite[Theorem~VI.1.33]{K}
the form $-\lambda\frd-\fra_0$ is a closed sectorial form with domain
$ H_{\frac{1}{2}}$ and hence the same is true for $-\frt(\lambda)$.
\end{proof}

\noindent
Since, by Theorem~\ref{rel_bd}, the forms $-\frt(\lambda)$ are closed and sectorial for $\lambda\in\Omega$,
it follows from \cite[Theorem~VI.2.1]{K} that there exist m-sectorial  operators $-T(\lambda)$ 
representing $-\frt(\lambda)$, i.e.\
\begin{equation}\label{repr_frt}
  \frt(\lambda)[x,y] = \langle T(\lambda)x,y\rangle, \qquad x\in\cD(T(\lambda)),\, y\in H_{\half}.
\end{equation}
The operator function $\lambda\mapsto-T(\lambda)$, defined on $\Omega$, is a 
\emph{holomorphic family of type}~(B) (see \cite[\S VII.4]{K}) because 
the domain $\cD(\frt(\lambda))$ of the quadratic form $\frt(\lambda)$ is independent of $\lambda$ 
and $\lambda\mapsto\frt(\lambda)[x]$ is holomorphic on $\Omega$ for every $x\in\cD(\frt(\lambda))$.

The following theorem collects some properties of the operator family $T(\lambda)$, $\lambda\in\Omega$.
Let us first recall some definitions.
The \emph{spectrum} of the operator function $T$ is defined as follows:
\[
  \sigma(T) \defeq \bigl\{\lambda\in\Omega \mid T(\lambda) \text{ is not bijective} \bigr\}
  = \bigl\{\lambda\in\Omega \mid 0\in\sigma(T(\lambda))\bigr\}.
\]
Similarly, the \emph{essential spectrum} of the operator function $T$ is defined as
\[
  \sess(T) \defeq \bigl\{\lambda\in\Omega \mid T(\lambda) \text{ is not Fredholm}\bigr\}
  = \bigl\{\lambda\in\Omega \mid 0\in\sess(T(\lambda))\bigr\}.
\]

A number $\lambda\in\Omega$ is called an \emph{eigenvalue} of the operator function $T$
if there exists an $x\in\cD(T(\lambda))$, $x\ne0$, such that $T(\lambda)x=0$.
The \emph{point spectrum} is the set of all eigenvalues:
\begin{align*}
  \sigma_{\rm p}(T) \defequ \bigl\{\lambda\in\Omega \mid
  \exists\,x\in\cD(T(\lambda)),\,x\ne0,\,T(\lambda)x=0\bigr\} \\[1ex]
  &= \bigl\{\lambda\in\Omega \mid 0\in\sigma_{\rm p}(T(\lambda))\bigr\},
\end{align*}
where $\sigma_{\rm p}(T(\lambda))$ denotes the point spectrum of the operator
$T(\lambda)$ for fixed $\lambda\in\Omega$.
The \emph{geometric multiplicity} of an eigenvalue $\lambda$ of the operator
function $T$ is defined as the dimension of $\ker T(\lambda)$.

\begin{thm}\label{Ilmenau2}
Let $-T(\lambda)$ be the m-sectorial operator in $H$ associated with the form
$-\frt(\lambda)$ satisfying \eqref{repr_frt}.  Then, for all $\lambda \in \Omega$, we have
\[
  \dom(T(\lambda)) =  \bigl\{ x \in H_{\frac12} \bigm| \lambda Dx + A_0x \in H \bigr\}
\]
and, for $x\in \dom(T(\lambda))$,
\[
  T(\lambda)x = \lambda^2 x + \lambda D x + A_0x.
\]
The family $-T(\lambda)$ of m-sectorial operators, $\lambda \in \Omega$,
is a holomorphic family of type $(${\rm B}$)$.
For $\lambda\in\Omega\cap\RR$ the operators $T(\lambda)$ are self-adjoint
and bounded above.

In addition, we have $\sess(T) =\emptyset$.
\end{thm}

\begin{proof}
The form of the domain and the action of $T(\lambda)$ follow from \cite[Theorem~VI.2.1]{K}.
As already mentioned above, $-T$ is a holomorphic family of type (B).
For $\lambda\in\Omega\cap\RR$ the form $\frt(\lambda)$ is symmetric and
hence $T(\lambda)$ is self-adjoint; $-T(\lambda)$ is bounded below since it is m-sectorial. 

It remains to show that $\sess(T)=\emptyset$. 
We show that $T(\lambda)$ is semi-Fredholm with $\dim(\ker T(\lambda))<\infty$
for every $\lambda\in\Omega$.  Suppose that this is not the case.
Then there exists a singular sequence $(x_n)_{n\in\NN}$, i.e.\  $x_n\in\mathcal D(T(\lambda))$, $\|x_n\|=1$, $x_n \rightharpoonup 0$ and
$T(\lambda)x_n\to0$ as $n\to \infty$; see, e.g.\ \cite[Theorem~IX.1.3]{EE}. Set
\begin{equation*}
  M \defeq \frac{a+b|\lambda|^2+1}{1-b}
\end{equation*}
Decompose $x_n$ as follows:
\[
  x_n = u_n + v_n, \qquad u_n\in \cL_{[0,M]}(A_0),\; v_n\perp \cL_{[0,M]}(A_0).
\]
The sequence $u_n$ converges weakly to $0$ because the projection onto $\cL_{[0,M]}(A_0)$ is
weakly continuous.  Since $\cL_{[0,M]}(A_0)$ is finite-dimensional, $(u_n)$ and $(A_0u_n)$ converge
also strongly in $H$ to $0$, which implies
with \eqref{a0_rel_bdd_d} that
\begin{align*}
  \bigl|\frt(\lambda)[x_n]\bigr| 
  &= \big|\lambda^2 + \lambda\frd[x_n] + \fra_0[x_n]\big|
  \ge |\lambda|\frd[x_n] - |\lambda|^2 - \fra_0[x_n]
  \\[1ex]
  &\ge \Bigl(\frac{1}{b}-1\Bigr)\fra_0[x_n]-\frac{a}{b}-|\lambda|^2
  \\[1ex]
  &= \Bigl(\frac{1-b}{b}\Bigr)\langle A_0u_n+A_0v_n,u_n+v_n\rangle_{H_{-\frac{1}{2}}\times H_{\frac{1}{2}}}
  - \Bigl(\frac{a}{b} + |\lambda|^2\Bigr)
  \\[1ex]
  &\ge \Bigl(\frac{1-b}{b}\Bigr) M\|v_n\|^2
  - \Bigl(\frac{a}{b} + |\lambda|^2\Bigr)\\[1ex]
  &= \Bigl(\frac{a}{b} +|\lambda|^2 \Bigr) \Bigl(\|v_n\|^2-1\Bigr)
  + \frac{1}{b} \|v_n\|^2.
\end{align*}
As $\|v_n\|^2 \to 1$ and $\frac{1}{b}>1$, there exists $N\in \mathbb N$ such that
\[
  |\frt(\lambda)[x_n]|\ge 1 \qquad \text{if} \;\;  n\ge N.
\]
This is a contradiction to the fact that $T(\lambda)x_n\to0$.
Hence $T(\lambda)$ is semi-Fredholm with $\dim(\ker T(\lambda))<\infty$ for every $\lambda\in\Omega$.
It follows from \cite[Theorem~VI.2.5]{K} that $(T(\lambda))^*=T(\overline\lambda)$
and hence $\dim(\ran T(\lambda))^\perp = \dim(\ker T(\overline\lambda)) < \infty$,
which implies that $T(\lambda)$ is a Fredholm operator for every $\lambda\in\Omega$.
\end{proof}

\section{The associated operator $\mathcal A$}\label{sec:matrix}

Let 
\[
  \cH \defeq H_{\frac{1}{2}} \times H
\]
equipped with the standard inner product which is denoted by $\langle \cdot, \cdot \rangle$.

Equation \eqref{sys} is equivalent to the following standard first-order equation
\begin{equation*}
  \dot{x}(t) = {\mathcal A} x(t), 
\end{equation*}
where $\cA: \dom(\cA)\subset\cH\rightarrow\cH$ is given by
\begin{align}
  & \cA = \begin{bmatrix} 0 & I \\[0.5ex] -A_0 & -D \end{bmatrix}, 
  \nonumber\\[1ex]
  & \dom(\cA) = \left\{\vect{x}{y} \in H_{\frac{1}{2}} \times H_{\frac{1}{2}} \Bigm|
  A_0x + Dy \in H  \right\}.
  \label{domA}
\end{align}

The operator $\cA$ itself is not self-adjoint in the Hilbert space $\cH$.
It is easy to see (e.g.\ \cite{TWII}) that $\cA$ has a bounded inverse in $\cH$ given by
\begin{equation*} 
  \cA^{-1} = \begin{bmatrix} -A_0^{-1}D & -A_0^{-1} \\[1ex] I & 0 \end{bmatrix},
\end{equation*}
where $A_0^{-1}D$ is considered as an operator acting in $H_{\frac{1}{2}}$.
This together with the fact that $J\cA$, where
\begin{equation*}
  J \defeq \begin{bmatrix} I & 0 \\[0.3ex] 0 & -I \end{bmatrix},
\end{equation*}
is a symmetric operator in the Hilbert space $\cH$, implies the
self-adjointness of $J\cA$ in $\cH$.
For $\smvect{x_1}{y_1},\smvect{x_2}{y_2} \in  H_{\frac{1}{2}}\times H$ we define
an indefinite inner product on $\cH$ by
\begin{equation} \label{Krein}
  \left[ \vect{x_1}{y_1},\vect{x_2}{y_2}\right] \defeq
  \left\langle J \vect{x_1}{y_1},\vect{x_2}{y_2}\right\rangle
  = \langle x_{1},x_{2} \rangle_{H_\frac{1}{2}}-\langle y_{1},y_{2} \rangle.
\end{equation}
Then $(\cH, \Skindef )$ is a Krein space and $\cA$ is a self-adjoint operator
with respect to $\Skindef$ (note that the latter is equivalent to the
self-adjointness of $J\cA$ in $\cH$). Hence $\sigma(\cA)$ is
symmetric with respect to $\mathbb R$; see, e.g.\ \cite[Theorem~VI.6.1]{B}.
For the basic theory of Krein spaces and operators acting therein we refer to \cite{AI,B,G22}.

In the following theorem we collect the above considerations and use \cite{JTW}.

\begin{thm}
\label{thm:cont}
If \textup{\textbf{(A1)--(A4)}} hold, then
the operator $\cA$ is self-adjoint in the Krein space
$(\cH, \Skindef )$, its spectrum is contained
in the closed left half plane and is symmetric with respect to the real line. 
The operator $\cA$ has a bounded inverse and it is the generator of an analytic semigroup of contractions on the space $\cH$.
\end{thm}

In the following proposition we collect consequences from standard Fredholm theory and \cite[Theorem~4.1 and Corollary~5.2]{jatr}.

\begin{prop} \label{Ilmenau}
Assume that \textup{\textbf{(A1)--(A4)}} hold. Consider $A_0^{-1}D$
as an operator acting in $H_{\frac{1}{2}}$ and $A_0^{-1/2}DA_0^{-1/2}$
as an operator acting in $H$.
The spectrum in $\CC\setminus\sess(\cA)$ is a discrete set consisting only of
eigenvalues, and the set $\sigma(\cA)\setminus\RR$ has no finite accumulation point.  Moreover,
\begin{equation*} 
  \min  \sess(\mathcal A) = -\frac{1}{\delta_{\rm ess}},
  \qquad
  \max  \sess(\mathcal A) = -\frac{1}{\gamma_{\rm ess}}
\end{equation*}
and 
\begin{align*}
  \sess(\mathcal A) &= \left\{\lambda\in\mathbb C\backslash\{0\}\,\Bigm|\,
    \frac{1}{\lambda}\in \sess(-A_0^{-1}D)\right\}\\[1ex]
  &= \left\{\lambda\in\mathbb C\backslash\{0\}\,\Bigm|\,
    \frac{1}{\lambda}\in \sess\bigl(-A_0^{-1/2}DA_0^{-1/2}\bigr)\right\}.
\end{align*}
\end{prop}
On $\Omega$ the spectral problems of $\mathcal A$ and $T$ coincide, this is the content of the following proposition. 
However, the proof of Proposition~\ref{Flughafen}
is identical to the one in \cite[Proposition~4.4]{JLT16}; only the domain of definition of the function $T$ has to be changed.

\begin{prop} \label{Flughafen}
Assume that {\rm\textbf{(A1)--(A4)}} are satisfied. 
On $\Omega$ the spectra of $\mathcal A$ and $T$ coincide, i.e.\
\begin{equation*}
  \sigma_{\rm p}(\mathcal A) \cap  \Omega =
  \sigma(\mathcal A) \cap  \Omega = \sigma(T)=
  \sigma_{\rm p}(T).
\end{equation*}
For  $\lambda_0\in \sigma_{\rm p}(\mathcal A) \cap \Omega$
the geometric multiplicities coincide,
\[
  \dim \ker (\mathcal A - \lambda_0) = \dim \ker T(\lambda_0).
\]
\end{prop}

In the following proposition we show that the operator  $\cA-c$ has a finite number of negative
squares  for each $c\in\Omega\cap\RR$. Recall, that a self-adjoint
operator $S$ in a Krein space $(\cK,\Skindef)$ is said to have \emph{$\kappa$ negative
squares} for $\kappa \in \mathbb N$ if the Hermitian form $[S\,\cdot\,,\cdot\,]$,
defined on  $\dom(S)$ has $\kappa$ negative squares, that is, there exists a subspace of
dimension $\kappa$ such that $[S\,\cdot\,,\cdot\,]$ restricted to this
subspace is negative definite and $\kappa$ is the largest number
with this property; see, e.g.\ \cite{L82}. In an obvious way,
this definition extends to sesquilinear forms. Correspondingly, one defines $\kappa$ positive squares.

\begin{prop}\label{kappanegsq}
Let $c\in\Omega\cap\RR$.  Then the operator $\cA-c$ has a finite number $\kappa_c$ of
negative squares in the Krein space  $H_{\half}\times H$ with the
inner product \eqref{Krein}, where the number of negative squares is given by
\begin{equation}\label{pos_squ}
  \kappa_c = \dim\cL_{(0,\infty)}\bigl(T(c)\bigr).
\end{equation}
This implies, that the non-real spectrum of $\mathcal A$ consists
of at most $\kappa$ pairs $\{\mu_i,\overline\mu_i\}$ of eigenvalues with finite-dimensional
algebraic eigen\-spaces, where $\kappa$ is the smallest of all
numbers in \eqref{pos_squ},
\begin{equation}\label{pos_squ2}
  \kappa \defeq \min_{c\in\Omega\cap\RR} \kappa_c.
\end{equation}
There are at most $\kappa$ different real  eigenvalues of $\mathcal A$ with
corresponding Jordan chains of length greater than one, each of them having length at most $2\kappa +1$.
\end{prop}

\begin{proof}
We start by showing that the quantity in \eqref{pos_squ} is a finite number.
Let $c\in\Omega\cap\RR=\bigl(-\infty,-\frac{1}{\delta_{\textup{ess}}}\bigr)$.
It follows from Theorem~\ref{rel_bd} that there exist constants $a,b$ with $a\ge0$ and $0<b<1$ 
such that \eqref{a0_rel_bdd_d} holds with $\lambda$ replaced by $c$. Hence, for $x\in H_{\frac{1}{2}}$, we have
\begin{align}
  \nonumber
  \frt(c)[x] &= c^2\|x\|^2-|c|\frd[x]+\fra_0[x]
  \\[0.5ex]
  \nonumber
  &\le c^2\|x\|^2 - \frac{1}{b}\fra_0[x] + \frac{a}{b}\|x\|^2 + \fra_0[x]
  \\[0.5ex]
  \label{baba}
  &= -\Bigl(\frac{1}{b}-1\Bigr)\fra_0[x]
  + \Bigl(c^2+\frac{a}{b}\Bigr)\|x\|^2.
\end{align}
Since the expression $\frac{1}{b}-1$ is positive and $A_0$ is non-negative
and has compact resolvent, the form $\frt(c)$ has a finite number of positive squares, 
which equals $\dim\cL_{(0,\infty)}\bigl(T(c)\bigr)$. 

Let $\kappa_c$ be given by \eqref{pos_squ}.  We show that $\cA-c$ has $\kappa_c$ negative squares.
For $\smvect{x}{y}\in\cD(\cA)$ we have
\begin{align} 
  & \frac{1}{c}\frt(c)[y]-c\Bigl\|x-\frac{1}{c}y\Bigr\|_{H_\half}^2 
  \nonumber\\[1ex]
  &= \frac{1}{c}\bigl\langle c^2y+cDy+A_0y,y\bigr\rangle_{H_{-\half}\times H_\half}
  - c\langle A_0x,x\rangle_{H_{-\half}\times H_\half}
  + \langle A_0x,y\rangle_{H_{-\half}\times H_\half} 
  \nonumber\\[1ex]
  &\quad + \langle A_0y,x\rangle_{H_{-\half}\times H_\half}
  -\frac{1}{c}\langle A_0y,y\rangle_{H_{-\half}\times H_\half} \notag\\[1ex]
  &= \bigl\langle A_0(-cx+y),x\bigr\rangle_{H_{-\half}\times H_\half}
  + \bigl\langle A_0x+Dy+cy,y\bigr\rangle_{H_{-\half}\times H_\half} 
  \nonumber\\[1ex]
  &= \bigl\langle -cx+y,x\bigr\rangle_{H_\half} + \bigl\langle A_0x+Dy+cy,y\bigr\rangle
  \nonumber\\[1ex]
  &= \left[\begin{bmatrix} -c & I \\[1ex] -A_0 & -D-c \end{bmatrix}
  \vect{x}{y}, \vect{x}{y}\right]  
  = \left[(\cA-c)\vect{x}{y},\vect{x}{y}\right],
  \label{peak}
\end{align}
where we have used that $A_0x+Dy\in H$ by \eqref{domA}.
Since $c<0$, the expression $-c\|x-c^{-1}y\|_{H_\half}^2$ is always non-negative.
Let $\cL \subset \dom (\cA)$ be a subspace such that
$[(\cA-c) \cdot , \cdot ]$ restricted to $\cL$ is  negative definite. Then the above calculation shows that
$\frt(c)[y]>0$ for all non-zero $y$ from the space
\[
  \widehat\cL\defeq\biggl\{y\in H: \exists\,x\text{ such that }
  \vect{x}{y}\in\cL\biggr\}.
\]
Obviously, the dimension of $\widehat \cL$ is at most the dimension of $\cL$.
If $\dim\widehat\cL<\dim \cL$, then there exists an $x\ne0$ such that
$\smvect{x}{0}\in\cL$, and thus
\[
  \biggl[\bigl(\cA-c\bigr)\vect{x}{0},\vect{x}{0}\biggr] = -c\|x\|_{H_\half}^2 > 0,
\]
which is a contradiction. Hence one has $\dim\widehat \cL=\dim \cL$.  
Since the form $\frt[\cdot]$ has $\kappa_c$ positive squares, the dimension of $\cL$ is bounded by $\kappa_c$. 
Therefore, the operator $\cA-c$ has at most $\kappa_c$ negative squares.

As $\cL_{(0,\infty)}\bigl(T(c)\bigr)$ is a $\kappa_c$-dimensional subspace, it is a subset of
$\cD\bigl(T(c)\bigr)$, i.e.\ by Theorem~\ref{Ilmenau2},
for each $y\in\cL_{(0,\infty)}\bigl(T(c)\bigr)$,
\begin{equation}\label{nott}
  cDy + A_0y \in H.
\end{equation}
Set
\[
  \wt\cL \defeq \left\{\vect{\frac{1}{c}y}{y}: y\in
  \cL_{(0,\infty)}\bigl(T(c)\bigr)\right\},
\]
which satisfies $\wt\cL\subset\cD(\cA)$ by \eqref{nott}.
Moreover, for $\smvect{x}{y}\in\wt\cL\setminus\{0\}$, it follows from \eqref{peak} that
\[
  \biggl[\bigl(\cA-c\bigr)\vect{x}{y},\vect{x}{y}\biggr]
  = \biggl[\bigl(\cA-c\bigr)\vect{\frac{1}{c}y}{y},\vect{\frac{1}{c}y}{y}\biggr]
  = \frac{1}{c}\frt(c)[y] < 0.
\]
Since $\dim\wt\cL=\kappa_c$, this shows that $\cA-c$ has at least $\kappa_c$ negative
squares.  Hence the number of negative squares of $\cA-c$ is equal to $\kappa_c$.

The remaining assertions of Proposition \ref{kappanegsq} follow from
\cite{L82} (or, for a more extended proof, we refer to \cite{AKT} and \cite[Theorem 3.1]{BT07}).
\end{proof}

\begin{rem}\label{rem_2}
We mention that, apart from the formula \eqref{pos_squ2} for the number $\kappa$,
the statements in Proposition~\ref{kappanegsq} are essentially
contained in \cite[Corollary~4.4]{JTW}, where they are proved differently.
Moreover, in  \cite{JTW} sufficient conditions are given for the
existence of a Riesz basis of $H_{\frac{1}{2}} \times H$ consisting of eigenvectors
 and finitely many associated vectors of $\cA$.
\end{rem}

\section{Variational Principles}\label{sec:varprinciple}

In the following we prove variational principles for real eigenvalues of $\mathcal A$.
To this end we introduce functions $p_+$ and $p_-$ which will serve as generalized Rayleigh functionals for $T$ on appropriate intervals.
For fixed $x\in H_{\frac{1}{2}}\backslash\{0\}$ consider the equation
\begin{equation}\label{scalar_qu_eq0}
  \frt(\lambda)[x] = \lambda^2\|x\|^2 + \lambda\frd[x] + \fra_0[x] = 0
\end{equation}
as an equation in $\lambda$.

\begin{defn}\label{Gustav2AAA}
If \eqref{scalar_qu_eq0} for $x\in H_{\frac12}\setminus\{0\}$
 does not have a real solution, then we set
\[
  p_+(x) \defeq -\infty, \qquad p_-(x) \defeq +\infty.
\]
Otherwise, we denote the solutions of \eqref{scalar_qu_eq0} by $p_\pm(x)$,
\begin{equation}\label{scalar_qu_eq17}
  p_\pm(x) \defeq \frac{-\frd[x] \pm \sqrt{\bigl(\frd[x]\bigr)^2 - 4\|x\|^2\fra_0[x]}\,}{2\|x\|^2}\,.
\end{equation}
Then the values of $p_+(x)$ and $p_-(x)$ belong to $(-\infty,0)\cup\{\pm\infty\}$. Set
\begin{equation*}
    \cD^* \defeq
 \bigl\{x\in H_{\frac12}\setminus\{0\} \mid p_\pm(x) \text{ are finite}\bigr\}.
\end{equation*}
\end{defn}

For the range of applicability of the variational principles in Theorem~\ref{theo111} below
we introduce the following numbers:
\begin{equation}\label{def_alpha_pm}
  \alpha_- \defeq \sup_{x\in\cD^*} p_-(x), \qquad
  \alpha_+ \defeq \inf_{x\in\cD^*} p_+(x).
\end{equation}

The following theorem is one of the main results of this paper.

\begin{thm}\label{theo111}
Assume that \textup{\textbf{(A1)--(A4)}} are satisfied.
Then $\cD^*\ne\emptyset$ and
\begin{equation}\label{upbdalpha}
  -\frac{2}{\delta} \le \alpha_+ \le -\frac{1}{\delta_{\textup{ess}}}, \qquad -\infty <\alpha_- \le -\frac{1}{\gamma}.
\end{equation}
Further, the following statements hold. 
\begin{itemize}
\item[\rm (i)]
Assume that $\alpha_-<-\tfrac{1}{\delta_{\rm ess}}$.
Then the set $\sigma(\cA)\cap \bigl(\alpha_-,-\tfrac{1}{\delta_{\rm ess}}\bigr)$ is either empty or consists of a finite or infinite
sequence of isolated semi-simple eigenvalues of finite multiplicity of $\cA$, which
can accumulate at most at $-\frac{1}{\delta_{\rm ess}}$.
Choose $c_+\in\bigl(\alpha_-,-\frac{1}{\delta_{\textup{ess}}}\bigr)$ 
and denote the eigenvalues in $\bigl[c_+,-\frac{1}{\delta_{\textup{ess}}}\bigr)$ by
$(\lambda_j^+)_{j=1}^{N^+}$, $N^+\in\NN\cup\{\infty\}$, in non-decreasing order,
counted according to their multiplicities: $\lambda_1^+\le \lambda_2^+\le \cdots$
if such eigenvalues exist, and set $N^+\defeq0$ otherwise.
Further, set
\[
  \kappa_{c_+} \defeq \dim\cL_{(0,\infty)}(T(c_+)).
\]
Then
\begin{equation}\label{minmax_lambda_pl}
\begin{split}
  \lambda_n^+ &= \min_{\substack{L\subset H_{1/2} \\ \dim L =\kappa_{c_+} +n}}
  \;\;\sup_{\substack{x\in L \\ x \ne 0}}\;\; p_+(x)\\[1ex]
  &= \max_{\substack{L\subset H \\ \dim L = \kappa_{c_+} +n-1}}
  \;\; \inf_{\substack{x\in H_{1/2} \\ x \ne 0,\, x \perp L}}
  \; p_+(x)
\end{split}
\end{equation}
for $n\in\NN$, $n\le N^+$.
\item[\rm (ii)]
The set $\sigma(\cA)\cap(-\infty,\alpha_+)$ consists of an infinite
sequence of isolated semi-simple eigenvalues of finite multiplicity, which accumulates only at $-\infty$.
Choose $c_-\in(-\infty,\alpha_+)$ and denote the eigenvalues in $(-\infty,c_-]$ 
by $(\lambda_j^-)_{j=1}^{\infty}$ in non-increasing order,
counted according to their multiplicities: $\lambda_1^-\ge \lambda_2^-\ge \cdots$.
Further, set
\[
  \kappa_{c_-} \defeq \dim\cL_{(0,\infty)}(T(c_-)).
\]
Then
\begin{equation}\label{minmax_lambda_mi}
\begin{split}
  \lambda_n^- &= \max_{\substack{L\subset H_{1/2} \\ \dim L = \kappa_{c_-}+n}}
  \;\; \inf_{\substack{x\in L \\ x \ne 0}}
  \;\; p_-(x) \\[1ex]
  &=  \min_{\substack{L\subset H \\ \dim L = \kappa_{c_-}+n-1}}
  \;\; \sup_{\substack{x\in H_{1/2} \\ x \ne 0,\, x \perp L}}
  \; p_-(x).
\end{split}
\end{equation}
for $n\in\NN$.
\end{itemize}
\end{thm}

\begin{figure}[ht]
\begin{center}
\begin{tikzpicture}[scale=1.0]
  \fill[red!18] (-11,-0.2) rectangle (-8,0.2);
  \fill[blue!14] (-6.5,-0.2) rectangle (0,0.2);
  \fill[red!30] (-8,0) rectangle (-6.5,0.2);
  \fill[blue!25] (-8,-0.2) rectangle (-6.5,0);
  \draw[->] (-11,0) -- (1,0);
  \draw (0,-0.2) -- (0,0.2);
  \node at (-0,-0.45) {$0$};
  \draw (-0.8,-0.1) -- (-0.8,0.1);
  \node at (-0.8,-0.45) {$-\gamma^{-1}$};
  \draw (-2,-0.1) -- (-2,0.1);
  \node at (-2,-0.45) {$-\gamma_{\textup{ess}}^{-1}$};
  \draw (-4,-0.1) -- (-4,0.1);
  \node at (-4,-0.45) {$-\delta_{\textup{ess}}^{-1}$};
  \draw (-6.5,-0.2) -- (-6.5,0.2);
  \node at (-6.5,-0.45) {$\alpha_-$};
  \draw (-8,-0.2) -- (-8,0.2);
  \node at (-8,-0.45) {$\alpha_+$};
  \draw (-6,-0.1) -- (-6,0.1);
  \node at (-6,-0.45) {$c_+$};
  \filldraw[black] (-5.5,0) circle (0.04);
  \node at (-5.5,0.45) {$\lambda_1^+$};
  \filldraw[black] (-4.8,0) circle (0.04);
  \node at (-4.8,0.45) {$\lambda_2^+$};
  \node at (-4.25,0.45) {$\cdots$};
  \draw (-8.7,-0.1) -- (-8.7,0.1);
  \node at (-8.7,-0.45) {$c_-$};
  \filldraw[black] (-9.2,0) circle (0.04);
  \node at (-9.2,0.45) {$\lambda_1^-$};
  \filldraw[black] (-10,0) circle (0.04);
  \node at (-10,0.45) {$\lambda_2^-$};
  \node at (-10.6,0.45) {$\cdots$};
  \filldraw[black] (-7,0) circle (0.04);
  \filldraw[black] (-7.5,0.6) circle (0.04);
  \filldraw[black] (-7.5,-0.6) circle (0.04);
  \filldraw[black] (-7.2,0.8) circle (0.04);
  \filldraw[black] (-7.2,-0.8) circle (0.04);
  \draw[line width=2, blue!80] (-4,0) -- (-2,0);
  \node at (-3,0.45) {$\sigma_{\textup{ess}}(\mathcal{A})$};
\end{tikzpicture}
\caption{\small The values of $p_+$ and $p_-$ are contained in the light blue and light red areas respectively.
In Theorem~\ref{theo111} one characterizes eigenvalues below $-\delta^{-1}_{\textup{ess}}$ 
and outside the interval $[\alpha_+, \alpha_-]$.
The thick black dots are the eigenvalues; there are at most finitely many non-real ones. The real ones accumulate at $-\infty$ and maybe at $-1/\delta_{\textup{ess}}$.}
\label{fig:spektrum}
\end{center}
\end{figure}
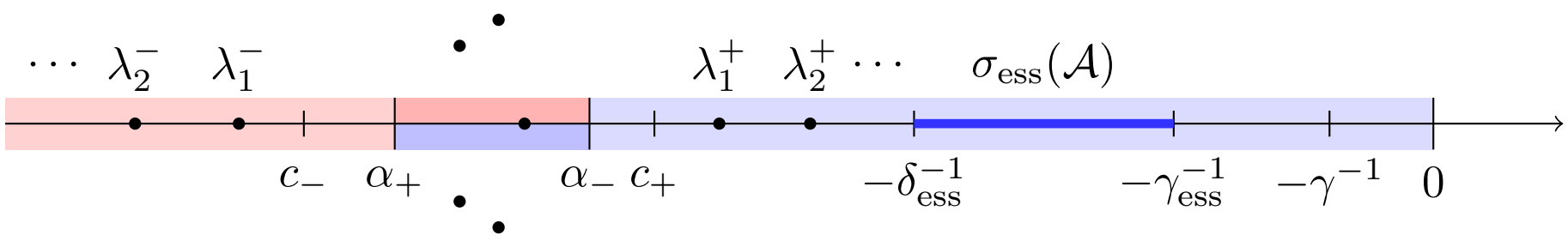

In Figure~\ref{fig:spektrum} we illustrate the location of $\alpha_\pm$, $c_\pm$, $-1/\delta_{\textup{ess}}$, 
$-1/\gamma_{\textup{ess}}$, $-1/\gamma$, the essential spectrum and some eigenvalues in the situation
when $\alpha_+<\alpha_-<-\frac{1}{\delta_{\textup{ess}}}$.

\begin{proof}
First we show that there exists a constant $M'>0$ such that, for $x\in H_{\half}\setminus\{0\}$,
\begin{equation}\label{a0geMp}
  \fra_0[x] \ge M'\|x\|^2 \quad\Rightarrow\quad \exists\,\lambda\in\Omega\cap\RR: \frt(\lambda)[x]=0.
\end{equation}
Fix $\lambda_0\in\Omega\cap\RR$.
By Theorem~\ref{rel_bd} there exist constants $a,b$ with $a\ge0$ and $0<b<1$ such that \eqref{a0_rel_bdd_d} holds 
(with $\lambda$ replaced by $\lambda_0$).  It follows from \eqref{baba} that, for $x\in H_{\half}\setminus\{0\}$
with $\fra_0[x]\ge M'\|x\|^2$,
\begin{equation*}
  \frt(\lambda_0)[x] 
  \le \Bigl(-\frac{1-b}{b}M'+\lambda_0^2+\frac{a}{b}\Bigr)\|x\|^2.
\end{equation*}
The right-hand side is less than zero if $M'$ is large enough.  As
$\lambda \mapsto\frt(\lambda)[x]$ describes a parabola,
there exists $\lambda\in(-\infty,\lambda_0]\subset\Omega\cap\RR$
such that $\frt(\lambda)[x]=0$.
This proves \eqref{a0geMp}, and, in particular, $\cD^*\ne\emptyset$ since $A_0$ is unbounded.

Next let us prove the inequalities in \eqref{upbdalpha}.
It follows from \eqref{nuevo} that, for $\lambda\in(-\frac{1}{\gamma},\infty)$ and $x\in H_\half\setminus\{0\}$,
\begin{equation*}
\begin{aligned}
  \frt(\lambda)[x] &= \lambda^2\|x\|^2 + \lambda\frd[x] + \fra_0[x] \\[1ex]
  &\ge \lambda^2\|x\|^2 - \frac{\frd[x]}{\gamma} + \fra_0[x] \ge \lambda^2\|x\|^2 \ge 0,
\end{aligned}
\end{equation*}
which implies that $p_\pm(x) \le -\frac{1}{\gamma}$ when $x\in\cD^*$.
Since $\cD^*\ne\emptyset$, this proves the second chain of inequalities in \eqref{upbdalpha}.

To prove $\alpha_+\ge -\frac{2}{\delta}$, let $x\in\cD^*$ with $\|x\|=1$.  
Then \eqref{nuevo2} implies
\[
  p_+(x) = \frac{\fra_0[x]}{p_-(x)} = -\frac{2\fra_0[x]}{\frd[x]+\sqrt{(\frd[x])^2-4\fra_0[x]}\,} 
  \ge -\frac{2\fra_0[x]}{\frd[x]} \ge -\frac{2}{\delta},
\]
which proves the first inequality in \eqref{upbdalpha}.

Next we prove $\alpha_+\le-\frac{1}{\delta_{\textup{ess}}}=\min\sess(\cA)$.
Let $\lambda\in\sess(\cA)$ and set
\[
  \cG \defeq \cL_{(0,2|\lambda|]}(A_0^{1/2}),
\]
which is a finite-dimensional subspace of $H_{\half}$.
Since $\lambda\in\sess(\cA)$, there exist $\binom{x_n}{y_n}$, $n\in\NN$,  
in $\dom(\cA) \cap (\cG\times\cG)^{\perp}$ with
\begin{equation}\label{Sbruecken}
  \bigg\|\binom{x_n}{y_n}\bigg\|_{H_{\frac{1}{2}} \times H} = 1
  \qquad \text{and} \qquad \lim_{n \to \infty}
  \bigg\|\bigl(\lambda I-\cA\bigr)\binom{x_n}{y_n}\bigg\|_{H_{\frac{1}{2}} \times H}
  = 0.
\end{equation}
The relation $x_n\in\cG^\perp$ implies that
\begin{equation}\label{inGperp}
  \|x_n\|_{H_\half} = \big\|A_0^{1/2}x_n\big\|
  \ge 2|\lambda|\,\|x_n\|.
\end{equation}
It follows from the second relation in \eqref{Sbruecken} that $\|y_n-\lambda x_n\|_{H_{\frac{1}{2}}}\to 0$ as $n\to \infty$.
Together with \eqref{inGperp} and the first relation in \eqref{Sbruecken}, we obtain
\begin{align}
  \big|\|y_n\|^2-\lambda^2\|x_n\|^2\big|
  &= \big|\|y_n\|-|\lambda|\|x_n\|\big|\cdot\bigl(\|y_n\|+|\lambda|\|x_n\|\bigr)
  \nonumber\\[0.5ex]
  &\le \|y_n-\lambda x_n\|\cdot\Bigl(\|y_n\|+\frac{1}{2}\|x_n\|_{H_\half}\Bigr)
  \nonumber\\[0.5ex]
  & \le \|y_n-\lambda x_n\|\cdot\Bigl( \|y_n\|+\frac{1}{2}
  \sqrt{1-\|y_n\|^2}\Bigr)
  \nonumber\\[0.5ex]
  &\le \frac{3}{2}\|y_n-\lambda x_n\| \to 0.
  \label{yn_laxn}
\end{align}
Again the first relation in \eqref{Sbruecken} and \eqref{inGperp} imply that
\begin{align*}
  \|x_n\|_{H_\half}^2 - \|y_n\|^2 &= \frac{3}{5}\Bigl(\|x_n\|_{H_\half}^2+\|y_n\|^2\Bigr)
  + \frac{2}{5}\Bigl(\|x_n\|_{H_\half}^2-4\|y_n\|^2\Bigr)
  \\[0.5ex]
  &\ge \frac{3}{5} + \frac{2}{5}\Bigl(4\lambda^2\|x_n\|^2-4\|y_n\|^2\Bigr),
\end{align*}
which, together with \eqref{yn_laxn} yields
\begin{equation}\label{SbrueckenII}
  \liminf_{n\to\infty}\Bigl(\|x_n\|_{H_{\frac{1}{2}}}^2 - \|y_n\|^2 \Bigr) \ge \frac{3}{5}.
\end{equation}
Further, $D:H_{\frac{1}{2}} \to H_{-\frac{1}{2}}$ is a bounded operator, 
which yields $D(y_n-\lambda x_n)\to 0$ in $H_{-\frac{1}{2}}$ as $n\to\infty$. 
The second relation in \eqref{Sbruecken} also implies that
\[
  (\lambda^2 + \lambda D + A_0)x_n 
  = -(D+\lambda)(y_n-\lambda x_n)+A_0x_n+(D+\lambda)y_n 
  \to 0,
\]
where the convergence is in $H_{-\frac12}$.
Since $x_n$ is bounded in $H_\half$ by \eqref{Sbruecken}, this implies that
\[
  \frt(\lambda)[x_n] = \bigl\langle(\lambda^2+\lambda D+A_0)x_n,x_n\bigr\rangle_{\mphhalf} \to 0.
\]
Let us set
\begin{align*}
  \varphi_n \defequ -\lambda x_n + y_n 
  \\[1ex]
  \psi_n \defequ -A_0 x_n - (D+\lambda)y_n,
\end{align*}
which, by \eqref{Sbruecken}, satisfy $\varphi_n\to 0$ in $H_\half$
and $\psi_n\to 0$ in $H$ as $n\to\infty$.
For the derivative of $\frt(\cdot)[x_n]$ we obtain
\begin{align*}
  & \frt'(\lambda)[x_n]
  = \bigl\langle \lambda x_n+(D+\lambda)x_n,x_n\bigr\rangle_{\mphhalf} \\[0.5ex]
  &= \langle y_n-\varphi_n,x_n\rangle
  + \frac{1}{\lambda}\bigl\langle(D+\lambda)(y_n-\varphi_n),x_n\bigr\rangle_{\mphhalf} \\[0.5ex]
  &= \frac{1}{\lambda}\langle y_n-\varphi_n,y_n-\varphi_n\rangle
  - \frac{1}{\lambda}\langle A_0x_n+\psi_n,x_n\rangle_{\mphhalf}
  - \frac{1}{\lambda}\bigl\langle(D+\lambda)\varphi_n,x_n\bigr\rangle_{\mphhalf} \\[0.5ex]
  &= \frac{1}{\lambda}\langle y_n,y_n\rangle - \frac{2}{\lambda}\Re\langle y_n,\varphi_n\rangle
  + \frac{1}{\lambda}\langle \varphi_n,\varphi_n\rangle - \frac{1}{\lambda}\langle A_0x_n,x_n\rangle_{\mphhalf}
  - \frac{1}{\lambda}\langle \psi_n,x_n\rangle \\[0.5ex]
  &\quad - \frac{1}{\lambda}\bigl\langle(D+\lambda)\varphi_n,x_n\bigr\rangle_{\mphhalf} \\[0.5ex]
  &= -\frac{1}{\lambda}\Bigl(\|x_n\|^2_{\hhalf}-\|y_n\|^2\Bigr) + \xi_n
\end{align*}
where
\begin{align*}
  \xi_n \defeq& -\frac{1}{\lambda}\biggl[2\Re\langle y_n,\varphi_n\rangle - \langle \varphi_n,\varphi_n\rangle
  + \langle \psi_n,x_n\rangle + \bigl\langle(D+\lambda)\varphi_n,x_n\bigr\rangle_{\mphhalf}\biggr] \\[1ex]
  &\to 0 \qquad\text{as}\;n\to\infty.
\end{align*}
From this and \eqref{SbrueckenII} it follows that
\[
  \liminf_{n\to\infty}\frt'(\lambda)[x_n] > 0.
\]
Relation \eqref{SbrueckenII} also implies that $\liminf_{n\to\infty}\|x_n\|_{\hhalf}>0$.  
Without loss of generality we can assume that $x_n\ne0$ and $\frt'(\lambda)[x_n]>0$ for all $n\in\NN$.
If we set
\[
  \tilde x_n \defeq \frac{1}{\|x_n\|_{\hhalf}}x_n,
\]
then there exists some $\eta>0$ such that
\begin{equation}\label{Berw01}
  \|\tilde x_n\|_{\hhalf} = 1, \quad \frt'(\lambda)[\tilde x_n] \ge \eta \quad\text{for all }n\in\NN
  \quad\text{and}\quad \lim_{n\to\infty}\frt(\lambda)[\tilde x_n] = 0.
\end{equation}
Set $\amin\defeq\min\sigma(A_0)$; then $\|\tilde x_n\|_{\hhalf}\ge\sqrt{\amin}\|\tilde x_n\|$
for all $n\in\NN$.
Now let $\eps\in(0,\frac{\eta\amin}{2})$ be arbitrary.  There exists an $n_0\in\NN$ such that
\begin{equation}\label{Berw02}
  \bigl|\frt(\lambda)[\tilde x_{n_0}]\bigr| \le \frac{\eps\eta}{2}.
\end{equation}
Using \eqref{Berw01} and \eqref{Berw02} we obtain
\begin{align*}
  \frt(\lambda+\eps)[\tilde x_{n_0}]
  &= \frt(\lambda)[\tilde x_{n_0}] + \eps\frt'(\lambda)[\tilde x_{n_0}] + \eps^2\|\tilde x_{n_0}\|^2 \\[0.5ex]
  &\ge -\frac{\eps\eta}{2} + \eps\eta > 0
\end{align*}
and
\begin{align*}
  \frt(\lambda-\eps)[\tilde x_{n_0}]
  &= \frt(\lambda)[\tilde x_{n_0}] - \eps\frt'(\lambda)[\tilde x_{n_0}] + \eps^2\|\tilde x_{n_0}\|^2 \\[0.5ex]
  &\le \frac{\eps\eta}{2} - \eps\eta + \frac{\eps^2}{\amin} = -\eps\Bigl(\frac{\eta}{2}-\frac{\eps}{\amin}\Bigr) < 0.
\end{align*}
By the Intermediate Value Theorem the function $\frt(\cdot)[\tilde x_{n_0}]$ has
a zero in the interval $(\lambda-\eps,\lambda+\eps)$ with positive derivative at that zero,
and hence $p_+(\tilde x_{n_0}) < \lambda+\eps$.
Since $\eps\in(0,\frac{\eta\amin}{2})$ was arbitrary, we have
$\alpha_+=\inf_{x\in\cD^*} p_+(x) \le \lambda$.
Now the claim follows because $\lambda$ was arbitrary in $\sess(\cA)$.

We now prove the statements in (i) and (ii).
Let us first show that $T$ satisfies assumption (II) in \cite{JLT16} with $\Omega$
as in \eqref{defOmegadelta0} and $\Delta=\Delta^+\defeq\bigl[c_+,-\frac{1}{\delta_{\textup{ess}}}\bigr)$ 
or $\Delta=\Delta^-\defeq(-\infty,c_-]$.  It follows from Theorem~\ref{Ilmenau2}
that $-T(\lambda)$ is a holomorphic family of type (B) and that $T(\lambda)$ is
self-adjoint for $\lambda\in\Delta^+$ or $\lambda\in\Delta^-$.
By Proposition \ref{kappanegsq} we have
\begin{equation*}
  \dim \mathcal L_{(0,\infty)}(T(\lambda))<\infty
\end{equation*}
for every $\lambda\in\Omega\cap\RR$, hence $\kappa_{c_+}$ and $\kappa_{c_-}$ are finite numbers.

Let us first consider the situation in (i).  On $\Delta^+$ the function
$\lambda\mapsto\frt(\lambda)[x]$ is increasing at value zero for $x\in H_\half\setminus\{0\}$
because this function is convex and if it has a zero
in $\Delta^+$, it is the greater one of the two zeros by the definition of $\alpha_-$
(note that a double-zero cannot lie in $\Delta^+$).
Therefore assumption $(\nearrow)$ in\cite{JLT16} with $\Delta=\Delta^+$
is fulfilled.
Moreover, the mapping $p_+$ is a generalized Rayleigh functional in the sense of \cite[Definition~2.1]{JLT16}.

By Proposition~\ref{Flughafen} the eigenvalues and their geometric multiplicities of $T$ and $\mathcal A$ 
coincide in $\Delta^+$.
The assertions apart from the semi-simplicity follow from \cite[Theorem~2.4]{JLT16}.
That the eigenvalues are semi-simple follows in the same way as in the proof of \cite[Theorem~4.8]{JLT16}.

The situation in (ii) is very similar. 
Any zero of $\lambda\to\frt(\lambda)[x]$ in $\Delta^-$
must be the smaller of the two zeros, and hence $(\searrow)$ in \cite{JLT16} is satisfied.
Moreover, the functional $p_-$ is a generalized Rayleigh functional in the sense of \cite[Definition~2.1]{JLT16}.
This shows that \cite[Theorem~2.6]{JLT16} can be applied, which yields most statements in (ii).
The semi-simplicity of the eigenvalues is shown in the same way as in the proof of \cite[Theorem~4.8]{JLT16}.

It remains to prove the accumulation of eigenvalues to $-\infty$.
Assume that $\cA$ has only a finite number $N^-$ of eigenvalues in $\Delta^-$.
Then, by \cite[(2.14)]{JLT16},
\[
  \sup_{\substack{L\subset H_{1/2} \\ \dim L = n}}
  \;\; \min_{x\in L\setminus\{0\}}
  \;\; p_-(x) = -\infty
\]
for $n>\kappa_{c_-}+N^-$.  This implies that, for all such $n$ and $n$-dimensional
subspaces $L$ of $H_{\frac{1}{2}}$, there exists an $x\in L\setminus\{0\}$ with
$p_-(x)=-\infty$, which is impossible by the definition of $p_-$, the implication in \eqref{a0geMp}
and the fact that $A_0$ is an unbounded operator.
Hence $N^-=\infty$ and $\lambda_n^-\to-\infty$ as $n\to\infty$.
\end{proof}

\bigskip

In the following theorem we consider the situation of a strongly damped or
hyperbolic operator polynomial, 
(see, e.g.\ \cite{M88}),
i.e.\ when $\frt(\lambda)[x]$ has always real
zeros.  

\begin{thm}\label{th:real_sp}
Suppose that \textup{\textbf{(A1)--(A3)}} are satisfied and assume
that, for every $x\in H_\half\setminus\{0\}$, the quadratic polynomial $\frt(\lambda)[x]$
has real zeros.
This is the case if and only if, for every $x\in H_\half\setminus\{0\}$,
we have
\begin{equation}\label{strong_damping_ml}
  \bigl\langle Dx,x\bigr\rangle_{\mphhalf}\geq 2\|x\|\|A_0^{1/2}x\|.
\end{equation}
Then $\cA$ satisfies also \textup{\textbf{(A4)}},
the spectrum of $\cA$ is real and there exists a number $\alpha_0\in(-\infty,-\frac{1}{\delta}]$ such that
\begin{equation}\label{hyperb1}
  \frt(\alpha_0)[x] \le 0, \qquad p_-(x) \le \alpha_0, \qquad p_+(x) \ge \alpha_0
\end{equation}
for all $x\in H_\half\setminus\{0\}$.  Let $\alpha_\pm$ be as in
\eqref{def_alpha_pm}.
Then
\begin{equation}\label{hyperb3}
  \alpha_- \le \alpha_0 \le \alpha_+ \le -\frac{1}{\deltaess}.
\end{equation}
In addition, the foolowing statements hold.
\begin{itemize}
    \item[\rm (i)] $\alpha_-<\alpha_+$ if and only if there exists some $\alpha_0'\in(-\infty,0)$
and some $\eta>0$ such that $\frt(\alpha_0')[x] \le -\eta\|x\|^2$ for all $x\in H_\half\setminus\{0\}$.
    \item[\rm (ii)] If $\alpha_-<\alpha_+$, then
    for any $c\in(\alpha_-,\alpha_+)$
\begin{equation*}
  (\alpha_-,\alpha_+)\subset\rho(\cA), \quad \alpha_+,\alpha_-\in\sigma(\cA),
  \quad\text{and}\quad \kappa_c = 0.
\end{equation*}
    \item[\rm (iii)] 
    If $\alpha_-=\alpha_+<-\frac{1}{\deltaess}$, then $\alpha_0=\alpha_-=\alpha_+$
is an eigenvalue of $\cA$.
\end{itemize}
\end{thm}

\begin{proof}
It follows from \cite[Lemma~2.1]{E72} that, if $\frt(\cdot)[x]$ has
real zeros for all $x\in H_\half\setminus\{0\}$, then there exists $\alpha_0\in(-\infty,0)$
such that
\begin{equation}\label{peterborough}
  p_-(x) \le \alpha_0 \le p_+(x) \qquad\text{for}\; x\in H_\half\setminus\{0\}.
\end{equation}
Since $\frt(\lambda)[x]$ is a quadratic polynomial in $\lambda$ with positive leading
coefficient, this implies that
\begin{equation}\label{hyperb10}
  \frt(\alpha_0)[x] = \alpha_0^2\|x\|^2 + \alpha_0\frd[x] + \fra_0[x] \le 0,
\end{equation}
which shows \eqref{hyperb1}.
It follows from \eqref{hyperb10} that, for $y\in H\setminus\{0\}$,
\begin{align*}
  \bigl\langle A_0^{-1/2}DA_0^{-1/2}y,y\bigr\rangle
  &= \bigl\langle DA_0^{-1/2}y,A_0^{-1/2}y\bigr\rangle_{\mphhalf}
  = \frd\bigl[A_0^{-1/2}y\bigr] \\[0.5ex]
  &\ge -\frac{1}{\alpha_0}\fra_0\bigl[A_0^{-1/2}y\bigr]
  = -\frac{1}{\alpha_0}\|y\|^2,
\end{align*}
which implies \textup{\textbf{(A4)}}. 
Relation \eqref{hyperb3} follows directly from \eqref{peterborough}, \eqref{upbdalpha}
and the definitions of $\alpha_+$ and $\alpha_-$.
The non-real spectrum of $\cA$ can consist only of eigenvalues. If $\lambda$
is a non-real eigenvalue of $\cA$ with eigenvector $\smvect{x}{y}$, then
$x\ne0$ and $\lambda$ is a non-real zero of $\frt(\cdot)[x]$, a contradiction.
Hence the spectrum of $\cA$ is real. 

If $\alpha_-<\alpha_+$, then, for $x\in H_{\frac12}$ with $\|x\|=1$
and $c\in(\alpha_-,\alpha_+)$,
\begin{align*}
  \frt(c)[x] = \bigl(c-p_+(x)\bigr)\bigl(c-p_-(x)\bigr)
  \le (c-\alpha_+)(c-\alpha_-) < 0,
\end{align*}
which shows that $c\in\rho(\cA)$ and $\kappa_c=0$ 
(see, e.g., \eqref{peak}). Hence, one can choose $\alpha_0'=c$.
Conversely, assume that $\frt(\alpha_0')[x]\le-\eta\|x\|^2$, $x\in H_{\frac12}$,
for some $\alpha_0'<0$ and $\eta>0$.  Then, for $x\in H_{\frac12}$ with $\|x\|=1$, we clearly
have $p_+(x)\ge\alpha_0'$.
Moreover, $\frt'(\alpha_0')[x]=2\alpha_0'+\langle Dx,x\rangle_{H_{-\frac12}\times H_{\frac12}}\ge2\alpha_0'$
and hence
\begin{align*}
  &\frt\Bigl(-\sqrt{\alpha_0'^2+\eta}\,\Bigr)[x] \\[1ex]
  &= \frt(\alpha_0')[x]+\frt'(\alpha_0')[x]\Bigl(-\sqrt{\alpha_0'^2+\eta}-\alpha_0'\Bigr)
  +\Bigl(-\sqrt{\alpha_0'^2+\eta}-\alpha_0'\Bigr)^2 \\[1ex]
  &\le -\eta+2\alpha_0'\Bigl(-\sqrt{\alpha_0'^2+\eta}-\alpha_0'\Bigr)
  +\Bigl(-\sqrt{\alpha_0'^2+\eta}-\alpha_0'\Bigr)^2 = 0,
\end{align*}
which implies that $p_-(x)\le-\sqrt{\alpha_0'^2+\eta}$.

Assume now that $\alpha_-<\alpha_+$ and let the notation be as in the proof of
Theorem~\ref{theo111}.  
If $\cA$ has at least one eigenvalue in the interval $(\alpha_-,-\frac{1}{\delta_{\rm ess}})$,
then
\[
  \lambda_1^+ = \min_{x\in \hhalf\setminus\{0\}} p_+(x) = \alpha_+
\]
by \eqref{minmax_lambda_pl}.  If $\cA$ has no eigenvalue in the
interval $(\alpha_-,-\frac{1}{\delta_{\rm ess}})$, then, by \cite[(2.8)]{JLT16} 
\[
  \inf_{x\in\hhalf\setminus\{0\}} p_+(x) \ge -\frac{1}{\delta_{\rm ess}}\,.
\]
Since the reverse inequality follows from \eqref{upbdalpha}, we have
equality and therefore $\alpha_+=\min\sess(\cA)$.
It follows from \eqref{minmax_lambda_mi} that
\[
  \lambda_1^- = \max_{x\in\hhalf\setminus\{0\}} p_-(x) = \alpha_-.
\]

Finally, assume that $\alpha_-=\alpha_+$.  Then $\alpha_0=\alpha_-=\alpha_+$
and hence $\frt(\alpha_0)[x]\le0$ for all $x\in H_{\hhalf}$, but there is
no $\eta>0$ such that $\frt(\alpha_0)[x]\le-\eta\|x\|^2$ for all $x\in H_{\hhalf}$ 
by what we have already proved.
This shows that $0\in\sigma(T(\alpha_0))$ and hence $\alpha_0\in\sigma(\cA)$.
\end{proof}

\begin{rem}
Let $\delta$ be as in \eqref{Elbersfeld}.
If
\begin{equation*}
  \delta^2\min\sigma(A_0) \ge 4,
\end{equation*}
then $\frt(\cdot)[x]$ has real zeros for each $x\in H_{\frac12}\setminus\{0\}$,
i.e.\ the assumption of Theorem~\ref{th:real_sp} is satisfied.
This follows because the expression under the square root in \eqref{scalar_qu_eq17},
the formula for $p_\pm(x)$, can be estimated with the help of \eqref{nuevo2}:
\[
  \bigl(\frd[x]\bigr)^2-4\|x\|^2\fra_0[x]
  \ge \delta^2\bigl(\fra_0[x]\bigr)^2-4\|x\|^2\fra_0[x]
  = \fra_0[x]\|x\|^2\Bigl(\delta^2\frac{\fra_0[x]}{\|x\|^2}-4\Bigr)
  \ge 0.
\]
\end{rem}

\section{Comparison of eigenvalues}\label{sec:comparison}

In the following theorem we use Theorem~\ref{theo111} to compare eigenvalues
of two operators.  Both an increase of $D$ and a decrease of $A_0$ moves the eigenvalues $\lambda_n^+$
to the right and the eigenvalues $\lambda_n^-$ to the left if enumerated appropriately.

\begin{thm}\label{comparison_strong}
Let operators $A_0$, $\hat A_0$, $D$ and $\hat D$ in the Hilbert space $H$ be given
so that $A_0$, $D$ and $\hat A_0$, $\hat D$, respectively, satisfy assumptions
\textup{\textbf{(A1)--(A4)}}, and let $\hat\cA$, $\hat T$, $\hat\frt$, $\hat\kappa_c$ 
be defined accordingly.
Assume that $H_{\frac12}=\hat H_{\frac12}$ and that, for $x\in H_{\frac12}$,
\begin{equation}\label{comp0_strong}
\begin{aligned}
  \langle\hat A_0x,x\rangle_{\hat H_{-\frac{1}{2}}\times \hat H_{\frac{1}{2}}}
  &\ge \langle A_0x,x\rangle_{H_{-\frac{1}{2}}\times H_{\frac{1}{2}}}, \\[1ex]
  \langle\hat Dx,x\rangle_{\hat H_{-\frac{1}{2}}\times \hat H_{\frac{1}{2}}}
  &\le \langle Dx,x\rangle_{H_{-\frac{1}{2}}\times H_{\frac{1}{2}}}.
\end{aligned}
\end{equation}
Then
\begin{equation}\label{comp001}
  \hat\delta_{\rm ess} \le \delta_{\rm ess}, \qquad
  \hat\alpha_+ \le \alpha_+, \qquad
  \hat\alpha_- \ge \alpha_-.
\end{equation}
\begin{enumerate}
\item[{\rm(i)}]
Assume that $\hat\alpha_-<-\tfrac{1}{\hat\delta_{\rm ess}}$.
Choose $c_+\in\bigl(\hat\alpha_-,-\frac{1}{\hat\delta_{\textup{ess}}}\bigr)$ 
and let $(\lambda_n^+)_{n=1}^{N^+}$ and $(\hat\lambda_n^+)_{n=1}^{\hat N^+}$,
$N^+, \hat N^+ \in\NN\cup\{\infty\}$, 
be the eigenvalues of $\cA$ and $\hat\cA$, respectively, in the interval $[c_+,-\frac{1}{\hat\delta_{\rm ess}})$,
both arranged in non-decreasing order, counted according to their multiplicities. 
If no such eigenvalues exists, set 
$N^+$ \textup{(}$\hat N^+$\textup{)} to zero.
Then $\hat\kappa_{c_+}\ge\kappa_{c_+}$ and
\[
  \hat\lambda_{n-\hat\kappa_{c_+}}^+ \le \lambda_{n-\kappa_{c_+}}^+
  \qquad \text{for }n\in\NN \ \text{with} \ n>\hat\kappa_{c_+},\,n\le\min\{N^+,\hat N^+\}.
\]
\item[{\rm(ii)}]
Let $c_-<\hat\alpha_+$ and let $(\lambda_n^-)_{n-1}^\infty$ and $(\hat\lambda_n^-)_{n=1}^\infty$
be the eigenvalues of $\cA$ and $\hat\cA$, respectively, in the interval $(-\infty,c_-]$,
both arranged in non-increasing order, counted according to their multiplicities. 
Then $\hat\kappa_{c_-}\ge\kappa_{c_-}$ and
\[
  \hat\lambda_{n-\hat\kappa_{c_-}}^- \ge \lambda_{n-\kappa_{c_-}}^-
  \qquad \text{for }n\in\NN \ \text{with} \ n>\hat\kappa_{c_-}.
\]
\end{enumerate}
\end{thm}

\begin{proof}
First we prove the inequality $\hat\delta_{\textup{ess}}\le\deltaess$.
Using the variational principle for self-adjoint operators we can characterise $\deltaess$ as follows:
\begin{align*}
  \deltaess &= \min\sess\bigl(A_0^{-1/2}DA_0^{-1/2}\bigr) \\[1ex]
  &= \sup_{n\in\NN}\;\; \inf_{\substack{L\subset H \\ \dim L=n}}\;\;
  \sup_{x\in L\setminus\{0\}}\;
  \frac{\bigl\langle A_0^{-1/2}DA_0^{-1/2}x,x\bigr\rangle}{\|x\|^2} \\[1ex]
  &= \sup_{n\in\NN}\;\; \inf_{\substack{L\subset H_{1/2} \\ \dim L=n}}\;\;
  \sup_{y\in L\setminus\{0\}}\;
  \frac{\langle Dy,y\rangle_{H_{-\frac{1}{2}}\times H_{\frac{1}{2}}}}{\langle
  A_0y,y\rangle_{H_{-\frac{1}{2}}\times H_{\frac{1}{2}}}}\,.
\end{align*}
Since a similar formula is true for $\hat\delta_{\textup{ess}}$, the inequalities in \eqref{comp0_strong} 
imply that $\hat\delta_{\textup{ess}}\le\deltaess$.

It also follows from the relations \eqref{comp0_strong} that, for negative $\lambda$,
\[
  \hat\frt(\lambda)[x] \ge \frt(\lambda)[x], \qquad x\in H_{\frac12},
\]
which implies that
\begin{equation}\label{comp002}
  \hat\kappa_c \ge \kappa_c \qquad\text{for all} \ c\in\Bigl(-\infty,-\frac{1}{\hat\delta_{\rm ess}}\Bigr)
\end{equation}
and
\begin{equation}\label{comp003}
  \hat p_-(x) \ge p_-(x) \qquad\text{and}\qquad
  \hat p_+(x) \le p_+(x), \qquad x\in\hhalf\setminus\{0\};
\end{equation}
note that $x\in\cD^*$ if $x\in\hat\cD^*$.
The remaining assertions follow from \eqref{comp002}, \eqref{comp003} and Theorem~\ref{theo111}.
\end{proof}

In the following example we consider the case when $D$ is a multiple of $A_0$.
This is used in Theorem~\ref	{comparison_with_example} below in combination with Theorem~\ref{comparison_strong}.

\begin{example}\label{ex:Dmult_of_A0}
Let $A_0$ be an operator in an infinite-dimensional space $H$ which
satisfies \textbf{(A1)} and \textbf{(A3)}.  Moreover, let $d>0$
and set 
\[
  D \defeq dA_0.
\]
Then also \textbf{(A2)} and \textbf{(A4)} are fulfilled
and $\gamma=\gamma_{\rm ess}=\delta=\delta_{\rm ess}=d$.
According to Proposition~\ref{Ilmenau} we have $\sess(\cA)=\bigl\{-\frac{1}{d}\bigr\}$.
A complex number $\lambda$ is an eigenvalue of $T$ if and only if there exists
a $\nu\in\sigma_{\rm p}(A_0)$ such that
\[
  \lambda^2+d\nu\lambda+\nu = 0.
\]
Hence, if $\nu_1\le\nu_2\le\cdots$ denote the eigenvalues of $A_0$, then
\begin{equation*}
  \sigma_{\rm p}(\cA) = \biggl\{-\frac{d\nu_n}{2}\pm\sqrt{\frac{d^2\nu_n^2}{4}-\nu_n}
  \;\bigg|\; n\in\NN\biggr\}.
\end{equation*}
It is easy to see that non-real eigenvalues (of which there are at most finitely many)
lie on the circle with centre at $(-\frac{1}{d},0)$ and radius $\frac{1}{d}$.
The real eigenvalues build two infinite sequences, one tending to $-\infty$
and the other one converging to $-\frac{1}{d}$, the only point in the essential spectrum.
If $\nu_1\ge\frac{4}{d^2}$, then \eqref{strong_damping_ml} is satisfied, 
all eigenvalues are real and
\begin{equation}\label{ex_alpha}
  \alpha_\pm = -\frac{d\nu_1}{2}\pm\sqrt{\frac{d^2\nu_1^2}{4}-\nu_1}.
\end{equation}
If $\nu_1<\frac{4}{d^2}$, then $\cA$ has non-real eigenvalues and, for $x\in\hhalf\setminus\{0\}$,
\[
  p_{\pm} = 
  \begin{cases}
    f_{\pm}\Bigl(\frac{\fra_0[x]}{\|x\|^2}\Bigr) & \text{if} \ \frac{\fra_0[x]}{\|x\|^2} \ge \frac{4}{d^2},
    \\[1ex]
    \mp\infty & \text{if} \ \frac{\fra_0[x]}{\|x\|^2} < \frac{4}{d^2},
  \end{cases}
\]
where
\[
  f_\pm(\nu) = -\frac{d\nu}{2} \pm \sqrt{\frac{d^2\nu^2}{4}-\nu}, \qquad 
  \nu\in\Bigl[\frac{4}{d^2},\infty\Bigr).
\]
Since $f_+$ is increasing and $f_-$ is decreasing on $\bigl[\frac{4}{d^2},\infty\bigr)$,
we have 
\[
  \alpha_\pm = f_\pm\Bigl(\frac{4}{d^2}\Bigr) = -\frac{2}{d}.
\]
\end{example}

\medskip

In the following main result of this section we establish explicit upper and lower bounds for the 
two eigenvalue branches by combining Theorem~\ref{comparison_strong} and Example~\ref{ex:Dmult_of_A0}
with $d=\delta$ and $d=\gamma$ respectively.

\begin{thm}\label{comparison_with_example}
Suppose that assumptions \textup{\textbf{(A1)--(A4)}} are satisfied and 
let $\nu_1\le\nu_2\le\cdots$ be the eigenvalues of $A_0$ in non-decreasing order.
Assume that 
$$
\delta^2\nu_1>4.
$$
Let $(\lambda_n^+)_{n=1}^{N_+}$ be the eigenvalues of $\cA$ in $\bigl(-\frac{\delta\nu_1}{2},-\frac{1}{\delta_{\textup{ess}}})$
in non-decreasing order,
let $(\lambda_n^-)_{n=1}^\infty$ be the eigenvalues of $\cA$ in $\bigl(-\infty,-\frac{\delta\nu_1}{2}\bigr)$
in non-increasing order. Then
\begin{alignat}{2}
  -\frac{\delta\nu_n}{2}+\sqrt{\frac{\delta^2\nu_n^2}{4}-\nu_n} &\le \lambda_n^+
  \le -\frac{\gamma\nu_n}{2}+\sqrt{\frac{\gamma^2\nu_n^2}{4}-\nu_n},
  \quad && n\in\{1,\ldots,N_+\},
  \label{est1}
  \\[1ex]
  -\frac{\gamma\nu_n}{2}-\sqrt{\frac{\gamma^2\nu_n^2}{4}-\nu_n} &\le \lambda_n^-
  \le -\frac{\delta\nu_n}{2}-\sqrt{\frac{\delta^2\nu_n^2}{4}-\nu_n},
  \quad && n\in\NN.
  \label{est2}
\end{alignat}
\end{thm}

\begin{proof}
We first apply Theorem~\ref{comparison_strong} with $\hat D=\delta A_0$ and $\hat A_0=A_0$.
It follows from \eqref{nuevo2} that the inequalities in \eqref{comp0_strong} are satisfied.
We obtain from \eqref{comp001} and \eqref{ex_alpha} that
\begin{align*}
  \alpha_+ \ge \hat\alpha_+ = -\frac{\delta\nu_1}{2}+\sqrt{\frac{\delta^2\nu_1^2}{4}-\nu_1} > -\frac{\delta\nu_1}{2},
  \\[1ex]
  \alpha_- \le \hat\alpha_- = -\frac{\delta\nu_1}{2}-\sqrt{\frac{\delta^2\nu_1^2}{4}-\nu_1} < -\frac{\delta\nu_1}{2}.
\end{align*}
Further, note that $\hat\delta_{\rm ess}=\delta$.
Let $N_+'$ be the number of eigenvalues of $\cA$ in $(-\frac{\delta\nu_1}{2},-\frac{1}{\delta})$.
Since $\hat\kappa_c=\kappa_c=0$ for $c\in(\hat\alpha_-,\hat\alpha_+)$, it follows from Theorem~\ref{comparison_strong}
that the second inequality in \eqref{est2} holds for $n\in\NN$ and the first inequality in \eqref{est1} 
holds for $n\in\{1,\ldots,N_+'\}$.
For $n\in\{N_+'+1,\ldots,N_+\}$ the first inequality in \eqref{est1} holds trivially because
\[
  -\frac{\delta\nu_n}{2}+\sqrt{\frac{\delta^2\nu_n^2}{4}-\nu_n} < -\frac{1}{\delta} \le \lambda_n^+.
\]
A second application of Theorem~\ref{comparison_strong} with $\hat D$ replaced by $D$ 
and $D$ replaced by $\gamma A_0$ yields the second inequality in \eqref{est1} and the first inequality in \eqref{est2}.
\end{proof}

\section{Application to an Euler--Bernoulli beam with Kelvin--Voigt damping}
\label{sec:beam}

We return to the Euler--Bernoulli beam introduced in the introduction in~\eqref{eq:beam}
with simply supported hinged boundary conditions. We assume that
\begin{equation}\label{eq:beam-coeff}
  C_d\in L^\infty(0,1),\qquad
  0 < c_-\defeq\essinf_{r\in(0,1)}C_d(r)
  \le c_+\defeq\esssup_{r\in(0,1)}C_d(r) < \infty.
\end{equation}
Choose $H=L^2(0,1)$ and let $A_0$ be the positive
self-adjoint realization of the fourth derivative with
\[
  \mathcal D(A_0)=\{u\in H^4(0,1):u(0)=u(1)=u''(0)=u''(1)=0\},
  \qquad A_0u=u''''.
\]
The  eigenvalues of $A_0$ are $(n\pi)^4$, $A_0^{-1}$ is compact and
\[
  H_{\frac{1}{2}} = \dom (A_0^{1/2}) = H^2(0,1)\cap H_0^1(0,1),  
  \qquad A_0^{1/2}u=-u''.
\]

The damping operator $D$ is a bounded operator
\[
  D:H_{\frac{1}{2}}\longrightarrow H_{-\frac{1}{2}}
\]
which is defined via the corresponding form,
\begin{equation}\label{Hymn}
  \langle Du,v\rangle_{H_{-\frac{1}{2}}\times H_\frac{1}{2}}
  \defeq
  \int_0^1 C_d(r)u''(r)\overline{v''(r)}\,dr,
  \qquad u,v\in H_\frac{1}{2}.
\end{equation}
Obviously, this form is symmetric and non-negative.
Regarding $D$ as an unbounded operator in \(L^2(0,1)\), its domain is
\[
  \cD(D)=\{u\in H_{\frac{1}{2}}  \mid
  C_du''\in H_{\frac{1}{2}}\},
\]
and \eqref{Hymn} gives
\[
  Du=A_0^{1/2} C_dA_0^{1/2} u,
  \qquad u\in \cD(D),
\]
or, what is the same,
\[
A_0^{-1/2}DA_0^{-1/2}= C_d.
\]
Therefore we have 
\begin{equation*}
 \sigma_{\rm ess}(A_0^{-1/2}DA_0^{-1/2})=\sigma(A_0^{-1/2}DA_0^{-1/2})=\essran C_d.
\end{equation*}
So \textbf{(A1)}--\textbf{(A4)}
are satisfied, with $\delta=\delta_{\rm ess}=c_-$ and $\gamma=\gamma_{\rm ess}=c_+$.
By Proposition~3.2, the essential spectrum
$\mathcal A$  is
\begin{equation*}
 \sigma_{\rm ess}(\mathcal A)
 =\{-1/s:s\in\essran C_d\}.
\end{equation*}
In particular, its left endpoint is $-1/c_-$.

\begin{cor}\label{cor:beam-bounds}
Assume \eqref{eq:beam-coeff} and
\begin{equation}\label{eq:beam-overdamped}
 c_-^2\pi^4>4.
\end{equation}
Let $N^+$ and the two sequences $\lambda_n^+$ and $\lambda_n^-$ be
numbered as in Theorem~\ref{comparison_with_example}, with $\delta=c_-$ and
$\nu_1=\pi^4$. Then, for $1\leq n\leq N^+$,
\begin{align*}
 -\frac{c_-(n\pi)^4}{2}
 +\sqrt{\frac{c_-^2(n\pi)^8}{4}-(n\pi)^4}
 &\leq \lambda_n^+\\
 &\leq
 -\frac{c_+(n\pi)^4}{2}
 +\sqrt{\frac{c_+^2(n\pi)^8}{4}-(n\pi)^4},\nonumber
\end{align*}
whereas, for every $n\in\mathbb N$,
\begin{align*}
 -\frac{c_+(n\pi)^4}{2}
 -\sqrt{\frac{c_+^2(n\pi)^8}{4}-(n\pi)^4}
 &\leq \lambda_n^-\\
 &\leq
 -\frac{c_-(n\pi)^4}{2}
 -\sqrt{\frac{c_-^2(n\pi)^8}{4}-(n\pi)^4}.\nonumber
\end{align*}
\end{cor}
\begin{proof}
The eigenvalues of $A_0$  and \eqref{eq:beam-overdamped}
satisfy all hypotheses of Theorem~\ref{comparison_with_example}. 
Substitution of $\delta=c_-$ and $\gamma=c_+$ into \eqref{est1}--\eqref{est2} gives the statement.
\end{proof}

\noindent
Birgit Jacob \\
School of Mathematics and Natural Sciences \\
University of Wuppertal \\
Gau\ss stra\ss e 20\\
D-42119 Wuppertal, Germany \\
Email: bjacob@uni-wuppertal.de

\bigskip

\noindent
Matthias Langer \\
Department of Mathematics and Statistics \\
University of Strathclyde \\
26 Richmond Street \\
Glasgow  G1 1XH, United Kingdom\\
Email: m.langer@strath.ac.uk

\bigskip

\noindent
Carsten Trunk \\
Institut f\"ur Mathematik \\
 Technische Universit\"at Ilmenau \\
Postfach 100565 \\
 D-98684 Ilmenau, Germany \\
 Email: carsten.trunk@tu-ilmenau.de


\begin{thebibliography}{99}

\bibitem{AKT}
R.\ Aydi, O.\ Kchaou, and C.\ Trunk,
Hermitian pencils and their representation
in Krein spaces, accepted in \textit{Integral Equations Operator Theory}.

\bibitem{AI}
T.\,Ya.~Azizov and I.\,S.~Iokhvidov,
\textit{Linear Operators in Spaces with an Indefinite Metric},
John Wiley \& Sons, 1989.

\bibitem{BSY}
H.\,T.~Banks, R.\,C.~Smith, and Y.~Wang,
The modeling of piezoceramic patch interactions with shells, plates, and beams,
\textit{Quart.\ Appl.\ Math.} \textbf{53} (1995), 353--381.

\bibitem{BT07}
J.~Behrndt and C.~Trunk,
On the negative squares of indefinite Sturm--Liouville operators,
\textit{J.\ Differential Equations} \textbf{238} (2007), 491--519.

\bibitem{B}
J.~Bognar,
\textit{Indefinite Inner Product Spaces},
Springer, 1974.

\bibitem{CLL}
S.~Chen, K.~Liu, and Z.~Liu,
Spectrum and stability for elastic systems with global or local Kelvin--Voigt damping,
\textit{SIAM J.\ Appl.\ Math.} \textbf{59} (1998), 651--668.

\bibitem{CourantHilbert}
R.~Courant and D.~Hilbert,
\textit{Methods of Mathematical Physics}, Wiley, 1953. 

\bibitem{Duffin1955}
R.\,J.~Duffin,
A minimax theory for overdamped networks,
\textit{J.\ Rational Mech.\ Anal.} \textbf{4} (1955), 221--233.

\bibitem{EE}
D.\,E.~Edmunds and W.\,D.~Evans,
\textit{Spectral Theory and Differential Operators},
Oxford University Press, 1987.

\bibitem{E72}
J.~Eisenfeld,
On symmetrization and roots of quadratic eigenvalue problems.
\textit{J.\ Funct.\ Anal.} \textbf{9} (1972), 410--422.


\bibitem{G22}
A. Gheondea, 
An Indefinite Excursion in Operator Theory: Geometric and Spectral Treks in Krein Spaces, 
Lond. Math. Soc. Lect. Note Ser. 476, Cambridge University Press, 2022.


\bibitem{JLT16}
B.~Jacob, M.~Langer, and C.~Trunk, 
Variational principles for self-adjoint operator functions arising from
second order systems, 
\textit{Oper.\ Matrices} \textbf{10} (2016), 501--531.


\bibitem{jatr}
B.~Jacob and C.~Trunk,
Location of the spectrum of operator matrices which are associated to
second order equations,
\textit{Oper.\ Matrices} \textbf{1} (2007), 45--60.

\bibitem{JT08} 
B.\ Jacob and C.\ Trunk, Spectrum and
analyticity of semigroups arising in elasticity theory and
hydromechanics, \textit{Semigroup Forum} \textbf{79} (2009), 79--100.


\bibitem{JTW}
B.\ Jacob, C.\ Trunk, and M.\ Winklmeier,
Analyticity and Riesz basis property of semigroups associated to damped vibrations,
\textit{J.\ Evol.\ Equ.} \textbf{8} (2008), 263--281.


\bibitem{K}
T.~Kato,
\textit{Perturbation Theory for Linear Operators},
Second Edition, Springer, 1976.

\bibitem{L82}
H.~Langer,
Spectral functions of definitizable operators in Krein spaces,
in: Functional Analysis (Dubrovnik, 1981),
\emph{Lecture Notes in Math.}, vol.~948, pp.~1--46, 1982.

\bibitem{M88} A.S.\ Markus, \textit{Introduction to the Spectral Theory of Operator Polynomials.} AMS Trans. Monographs, Providence, RI, 1988.


\bibitem{WT}
M.~Tucsnak and G.~Weiss,
How to get a conservative well-posed system out of thin air, Part I,
\textit{ESAIM Control Optim.\ Calc.\ Var.} \textbf{9} (2003), 247--274.

\bibitem{TWII}
M.~Tucsnak and G.~Weiss,
How to get a conservative well-posed system out of thin air, Part II,
\textit{SIAM J.\ Control Optim.} \textbf{42} (2003), 907--935.


\bibitem{Pu}
P.\,H.~You,
Boundary feedback control of elastic beam equation with structural damping and stability,
\textit{Acta Math.\ Appl.\ Sinica (English Ser.)} \textbf{6} (1990), 373--382.
\end{thebibliography}
\end{document}